\documentclass[12pt]{article}

\usepackage[fleqn]{amsmath}
\usepackage{amsfonts}
\usepackage{amssymb}
\usepackage{epsfig}
\usepackage{geometry}
\usepackage{ctable}
\usepackage{mathtools}

\usepackage[utf8]{inputenc}
\usepackage[T1]{fontenc}

\DeclarePairedDelimiter{\ceil}{\lceil}{\rceil}
\DeclarePairedDelimiter{\floor}{\lfloor}{\rfloor}

\newcommand{\Cset}{\mathbb{C}}

\newcommand{\Nset}{\mathbb{N}}

\newcommand{\Rset}{\mathbb{R}}

\newcommand{\cC}{\mathcal{C}}
\newcommand{\cS}{\mathcal{S}}

\newtheorem{theo}{Theorem}[section]

\newtheorem{coro}[theo]{Corollary}

\newtheorem{rema}[theo]{Remark}

\newcommand{\Cl}{{\mbox{P}_{\mbox{\tiny\sc{Cl}}}}}
\newcommand{\Clg}{{\mbox{Cl}_{p;s}^{}}}
\newcommand{\Clp}{{\mbox{Cl}_{\hfrac13;2}^{}}}

\newcommand{\sign}{{\mbox{sign}}}

\newcommand{\hfrac}[2]{{{#1}/{#2}}}

\begin{document}
\title{Order and Chaos in some  Deterministic \\ Infinite Trigonometric Products}
\author{\sc{Leif Albert \&\ Michael K.-H. Kiessling}\\ 
{\small{Department of Mathematics, Rutgers University, Piscataway, NJ 08854}}}
\date{Version of April 28, 2017} 
\maketitle

\begin{abstract}
\noindent
 It is shown that the deterministic infinite trigonometric products
\begin{equation}\notag
\prod_{n\in\Nset}\left[1- p +p\cos\left(\textstyle n^{-s}_{_{}}t\right)\right] =: \Clg(t) 
\end{equation}

\vspace{-10pt}\noindent
with parameters $p\in (0,1]\ \&\  s>\frac12$, and variable $t\in\Rset$, are inverse Fourier transforms of the 
probability distributions for certain random series $\Omega_{p}^\zeta(s)$ taking values in the real $\omega$ line; 
i.e. the $\Clg(t)$ are characteristic functions of the $\Omega_{p}^\zeta(s)$.
 The special case $p=1=s$ yields the familiar random harmonic series, while in general $\Omega_{p}^\zeta(s)$
is a ``random Riemann-$\zeta$ function,'' a notion which will be explained and illustrated --- and connected to the Riemann hypothesis.
 It will be shown that $\Omega_{p}^\zeta(s)$ is a very regular random variable,
having a 
probability density function (PDF) on the $\omega$ line which is a Schwartz function.
 More precisely, an elementary proof is given that 
there exists some $K_{p;s}^{}>0$, and 
a function $F_{p;s}^{}(|t|)$ bounded by $|F_{p;s}^{}(|t|)|\!\leq\! \exp\big(K_{p;s}^{} |t|^{1/(s+1)})$, 
and $C_{p;s}^{}\!:=\!-\frac1s\int_0^\infty\ln|{1-p+p\cos\xi}|\frac{1}{\xi^{1+1/s}}{\rm{d}}\xi$, such that
\begin{equation}\notag
\forall\,t\in\Rset:\quad 
\Clg(t) =
\exp\bigl({- C_{p;s}^{} \,|t|^{1/s}\bigr)F_{p;s}^{}(|t|)};
\end{equation}
the regularity of $\Omega_{p}^\zeta(s)$ follows.
 Incidentally, this theorem confirms a surmise by Benoit Cloitre, that 
$\ln \Clp(t) \sim -C\sqrt{t}\; \left(t\rightarrow\infty\right)$ for \emph{some} $C>0$.
 Graphical evidence suggests that $\Clp(t)$ is an empirically unpredictable (chaotic) function of $t$.
 This is reflected in the rich structure of the pertinent PDF (the Fourier transform 
of $\Clp$), and illustrated by random sampling of the Riemann-$\zeta$ walks, whose branching rules 
allow the build-up of fractal-like structures.

\end{abstract}
\vfill

\hrule\smallskip
\copyright{2017} The authors. Reproduction for non-commercial purposes is permitted.
\newpage

\section{Introduction and Summary}

 The Riemann hypothesis is perhaps the best-known open problem of mathematics. 
 It hypothesizes that all non-real zeros of Riemann's zeta function $\zeta(s)$, $s\in\Cset$, 
lie on the straight line $\frac12+i\Rset$, where $\zeta(s)$ is obtained from Euler's real (Dirichlet-)series 

\vspace{-20pt}
\begin{equation}\label{zetaFCT}
\zeta(s) = {\textstyle\sum\limits_{n\in\Nset}^{}}\; \frac{1}{n^s}, \quad s > 1,
\end{equation}
\vspace{-15pt}

\noindent
by analytic continuation to the complex plane; see \cite{Edwards} for a good introduction.
 The importance of Riemann's hypothesis derives from the fact that its truth would confirm deep putative insights into 
the distribution of the natural prime numbers --- a holy grail of number theory.
 This feat would also have applications: chiefly in encryption, 
but also in physics, see \cite{NHFG,SchuHu,Watkin}.
 It continues to fascinate the minds of professionals and amateurs alike.
 
 The latter group includes Benoit Cloitre, who has been documenting his experimental mathematical approach to 
number theory in general, and to the Riemann hypothesis in particular, on his homepage \cite{Cloitre}.
 Some years ago he pondered (``for no particular reason'')\footnote{\noindent{Private 
   communication by B.C. on 02.2016; we took the liberty to attach $_{\mbox{\tiny{\sc{Cl}}}}$ at Cloitre's~P$(t)$.}}
the deterministic infinite trigonometric product 

\vspace{-20pt}
\begin{equation}\label{CloitrePROD}
\textstyle\prod\limits_{n\in\Nset}\left[\frac23+\frac13\cos\left({t}/{n^{2}}\right)\right] =: \Cl(t),\quad t\in\Rset,
\end{equation}
\vspace{-15pt}

\noindent
which appears to be fluctuating chaotically about some monotone trend; see Fig.~1 and Fig.~2 below.
 Cloitre ``guessed'' that $\ln \Cl(t) \sim -C \sqrt{t}$ when $t\to\infty$ for \emph{some} constant $C>0$, 
which captures the trend asymptotically, and he asked us whether we can prove this. 
 The proof requires only elementary undergraduate mathematics and will be given in section \ref{THM}
(in fact, we prove a stronger result).
 But why does $\Cl(t)$ fluctuate apparently chaotically about its monotone trend?
 And what does this have to do with the Riemann hypothesis?
 Statistical physics offers some answers.

 We note (see section \ref{PROB}) that any trigonometric product

\vspace{-20pt}
\begin{equation}\label{genCLOITREprod}
\textstyle\prod\limits_{n\in\Nset}\left[1- p +p\cos\left({t}/{n^s}\right)\right] =: \Clg(t) ,\quad t\in\Rset, \quad p\in (0,1]\ \&\  s>\frac12, 
\end{equation}
\vspace{-15pt}

\noindent
is \emph{the characteristic function of a} ``\emph{random Riemann-$\zeta$ function}'' $\Omega^\zeta_{p}(s)$, i.e.
$\Clg(t)\equiv \Phi^{}_{\Omega^\zeta_p(s)}(t)$, where $\Phi^{}_{\Omega}(t) := \mbox{Exp}\big(e^{it\Omega}\big)$
with ``Exp'' denoting \emph{expected value}.
 Here,

\vspace{-20pt}
\begin{equation}\label{RANDOMzetaFCT}
\Omega^\zeta_{p}(s) 
:= {\textstyle\sum\limits_{n\in\Nset}^{}} R_p^{}(n)\frac{1}{n^s}, \quad s >\tfrac12,\quad p\in(0,1],
\end{equation}
\vspace{-15pt}

\noindent
where $\{R_p^{}(n)\in\{-1,0,1\}\}^{}_{n\in\Nset}$ is a sequence of i.i.d.
random coefficients, with Prob$(R_p^{}(n)=0) = 1-p$ and Prob$(R_p^{}(n)=1) = p/2=$ Prob$(R_p^{}(n)=-1)$.
 We draw heavily on the probabilistically themed publications by Kac 
\cite{Kac}, Morrison \cite{Morrison}, and Schmuland \cite{Schmuland}, 
in which \emph{the random harmonic series} ${\Omega^{\mbox{\tiny{harm}}}}\equiv\Omega^\zeta_{1}(1)$ is explored; in \cite{Schmuland} 
also the special case $\Omega^\zeta_{1}(2)$ is explored.
  We register that $p=1$ and $s=1$ in $\Clg(t)$ yields (cf. sect.5.2 in \cite{Morrison})

\vspace{-20pt}
\begin{equation}\label{MorrisonPRODharmonic}
\Phi^{}_{\Omega^{\mbox{\tiny{harm}}}}(t) = \textstyle\prod\limits_{n\in\Nset}^{}\cos \frac{t}{n},
\end{equation}

\vspace{-12pt}\noindent
while Cloitre's $\Cl(t)$ is the special case $p=\frac13$ and $s=2$ in $\Clg(t)$.

 Both $\zeta(s)$ and $-\zeta(s)$ are possible outcomes for such random Riemann-$\zeta$ functions $\Omega^\zeta_{p}(s)$, 
namely the extreme cases in which each $R_p^{}(n), n\in\Nset$, comes out $1$, respectively $-1$.
 While this is trivial, we anticipate that also $1/\zeta(s)$ is a possible outcome for $\Omega^\zeta_{p}(s)$, which is nontrivial
and going to be interesting!

 After introducing the notion of \emph{typicality} for the random walks associated to
$\Omega^\zeta_{p}(s)$ we will ask how typical $\zeta(s)$ and $1/\zeta(s)$ are.
 It should come as no surprise that $\zeta(s)$ is an extremely atypical outcome of a random Riemann-$\zeta$ walk, and so is $-\zeta(s)$.
 However, for the particular value of $p=6/\pi^2$, its reciprocal $1/\zeta(s)$ does exhibit several aspects of typicality.
 Intriguingly, as pointed out to us by Alex Kontorovich, \emph{if $1/\zeta(s)$ also exhibits a certain {particular} aspect 
of typicality, then the Riemann hypothesis is true, and false if not}! 
 This will be extracted from \cite{Edwards} in section \ref{typical}.

 Which of the many aspects of typicality are exhibited by $1/\zeta(s)$ is an interesting open question 
which may go beyond settling the Riemann hypothesis.
 We will use a paradox to argue, though, that $1/\zeta(s)$ cannot possibly exhibit each and every aspect of typicality, 
i.e. $1/\zeta(s)$ cannot be a perfectly typical random Riemann-$\zeta$ walk. 

 So much for the connection between Cloitre's $\Clp(t)$ and Riemann's $\zeta$ function. 
 
 As for our inquiry into Cloitre's surmise that $\ln \Clp(t) \sim -C\sqrt{t}\; \left(t\rightarrow\infty\right)$ for {some} $C>0$,
curiously some well-known probability laws emerged unexpectedly. 
 Using elementary analysis we will prove in section \ref{THM} that if 
$p\in(0,1]\ \&\ s> \frac12$, then there exists $K_{p;s}^{}\!>\!0$, and 
$F_{p;s}^{}(|t|)$ bounded by $|F_{p;s}^{}(|t|)| \leq \exp\big(K_{p;s}^{} |t|^{1/(s+1)}\big)$, so that

\vspace{-20pt}
\begin{equation}\label{genCLOITREtrendANDfluc}
\forall\,t\in\Rset:\quad 
\Clg(t) =
\exp\bigl(- C_{p;s}^{} \,|t|^{1/s}\bigr) F_{p;s}^{}(|t|)
\end{equation}

\vspace{-5pt}
\noindent
with \vspace{-10pt}
\begin{equation}\label{generalTRENDcoeff}
C_{p;s}^{} = -\frac1s\int_0^\infty\ln |{1-p+p\cos\xi}|\frac{1}{\xi^{1+1/s}}{\rm{d}}\xi;
\end{equation}

\vspace{-10pt}
\noindent
when $p\in(0,\frac12)$ the integral can be evaluated in terms of a rapidly converging series expansion.
 This result not only vindicates Cloitre's surmise as a corollary, we note
that the factor $\exp\bigl({- C_{p;s}^{} \,|t|^{1/s}}\bigr)$ at r.h.s.(\ref{genCLOITREtrendANDfluc}) in 
itself is a characteristic function --- of Paul L\'evy's \emph{stable laws}; see \cite{probBOOK}.
 Stable laws exist for all $s\geq 1/2$, but here $s = 1/2$ is ruled out because $C_{p;1/2}=\infty$.
 Be that as it may, stable L\'evy laws (which have applications in physics \cite{GaroniFrankel})
were discovered by answering a completely different question \cite{probBOOK,GaroniFrankel},
and the probabilistic reason why they would feature in the analysis of the random Riemann-$\zeta$ functions is presently obscure.

 Lest the reader gets the wrong impression that random Riemann-$\zeta$ functions were small perturbations of L\'evy random variables,
we emphasize that they are not! 
 Although the 
``chaotic factor'' $F_{p;s}^{}(|t|)$ in (\ref{genCLOITREtrendANDfluc}) is \emph{overwhelmed} 
by $\exp(-C_{p;s}^{}|t|^{1/s})$ when $|t|$ is large enough, 
$F_{p;s}^{}(|t|)$ is not approaching 1 and in fact responsible for relatively
large chaotic fluctuations of $\Clg(t)$ about the L\'evy trend; see Fig.~1\ \&\ Fig.~2.

\includegraphics[scale=0.44]{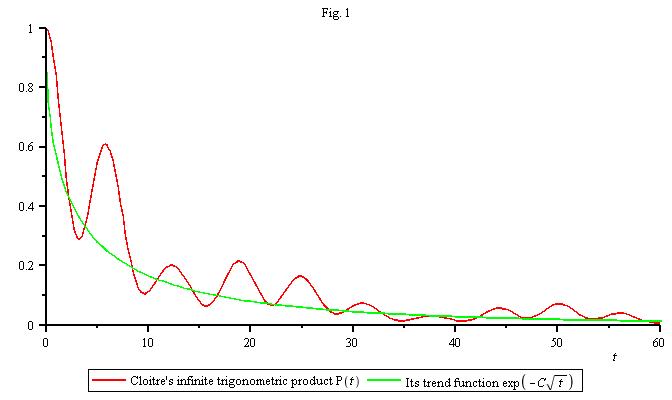}

\hspace{-20pt}
\includegraphics[scale=0.44]{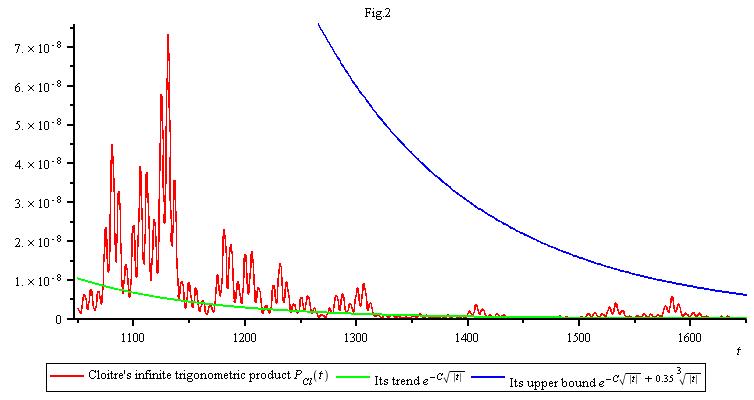}

 In section \ref{TandF} we will see that the ``empirically unpredictable'' behavior of $\Clp(t)$ is reflected in a 
\emph{fractal-like} structured probability distribution $\varrho^{\zeta}_{1/3;s}(d\omega)$ of $\Omega^\zeta_{1/3}(s)$ obtained by
Fourier transform  of $\Clp(t)$ (section \ref{PROB}).
 This is also illustrated in section \ref{RRZf} by random sampling of the Riemann-$\zeta$ walks.
 We will show, though, that $\varrho^{\zeta}_{p;s}(d\omega)$ is not supported on a true fractal. 
 Random variables supported on a fractal are discussed in \cite{DFT}, \cite{Morrison}, and \cite{PerezSchlagSolomyak};
see our Appendix on \emph{power walks}.

 The remainder of our paper supplies the details of our inquiry, and we conclude with a list of open questions.

\newpage

\section{Random Riemann-$\zeta$ functions}\label{RRZf}

 The random Riemann-$\zeta$ functions $\Omega_{p}^\zeta(s)$ given in (\ref{RANDOMzetaFCT}) have
random coefficients $R_p^{}(n)\in\{-1,0,1\}$ that can be generated by a two-coin tossing process.
 In this vein, let's write $R_p^{}(n) = \sigma(n)|R_p^{}(n)|$, where $\sigma(n)\in\{-1,1\}$ and $|R_p^{}(n)|\in\{0,1\}$. 
 One now repeatedly tosses both, a generally loaded coin with Prob$(H)=p\in(0,1]$ (where ``$H$'' means ``head''),
and a fair one, independently of each other and of all the previous tosses.
 The $n$-th toss of the generally loaded coin decides whether $|R_p^{}(n)|=0$ or $|R_p^{}(n)|=1$;
let's stipulate that $|R_p^{}(n)|=1$ when $H$ shows, which happens with probability $p$, and $|R_p^{}(n)|=0$ else.
 The concurrent and independent toss of the fair coin decides whether $\sigma(n)=+1$ or $\sigma(n)=-1$, either outcome being
equally likely.
 Incidentally, we remark that
the $R_{1/3}^{}(n)$ can also be generated by rolling a fair die --- if the $n$-th roll shows 1, then $R_{1/3}^{}(n)=1$, 
if it shows 6 then $R_{1/3}^{}(n)=-1$, and $R_{1/3}^{}(n)=0$ otherwise (which is the case $2/3$ of the time, in the long run).
 Also, it is clear that when $p=1$ then the loaded coin is superfluous, i.e. $R_1^{}(n)\in\{-1,1\}$ is generated with a single, 
fair coin.

 This completes the explanation of the ``experimental generation'' of our random Riemann $\zeta$ functions.
 Now let us understand which type of objects we have defined.\vspace{-10pt}

\subsection{Random Riemann-$\zeta$ walks and their kin}\label{RRZw}

 Evaluating a random Riemann-$\zeta$ function $\Omega_{p}^\zeta(s)$ for given $p\in(0,1]$ at any particular $s>\frac12$ 
turns (\ref{RANDOMzetaFCT}) into a numerical random series. 
 Recalling that an infinite series is defined as the sequence of its partial sums, viz.

\vspace{-18pt}
\begin{equation}\label{RZwDEF}
 \Omega_{p}^\zeta(s) = \left\{{\textstyle\sum\limits_{n=1}^N} R_p^{}(n)\frac{1}{n^s}\right\}^{}_{N\in\Nset}, 
\end{equation}
\vspace{-10pt}

\noindent
and interpreting $\sum_{n=1}^N R_p^{}(n)\frac{1}{n^s}$ as the position of a walker after $N$ random steps 
$R_p^{}(n)\frac{1}{n^s}, n=1,...,N$, we can identify such an evaluation of $\Omega_{p}^\zeta(s)$ for given $p\in(0,1]$ 
at a particular $s>\frac12$ with a \emph{random walk on the real $\omega$ line}.
 If the $n$-th toss of the pair of coins comes out on ``move,'' the walker moves $1/n^s$ units in the direction determined 
by the fair coin; otherwise he stays put (note that such a ``non-move'' is called a ``step,'' too).
 Starting at $\omega=0$, he keeps carrying out these random steps ad infinitum.
 We call this a ``random Riemann-$\zeta$ walk,'' and its outcome (whenever it converges)
is a ``random Riemann-$\zeta$ function'' evaluated at $s$.
 Absolute convergence is guaranteed for each and every such walk when $s>1$ (because the series (\ref{zetaFCT})
for $\zeta(s)$ converges absolutely for $s>1$), and by a famous result of Rademacher
conditional convergence holds with probability 1 when $s>\frac12$, see \cite{Kac}, \cite{Morrison}, and  \cite{Schmuland}.
 Since the harmonic series diverges logarithmically, the outcome of the random walks with $\frac12 <s\leq 1$
is distributed over the whole real line; see \cite{Schmuland} for $s=1$.

 To have some illustrative examples, we first pick $s=2$ and $p=\frac13$.
 In Fig.~3 we display (in black) the fractal tree (cf. \cite{Mandelbrot}, chpt.16; note its self-similarity) 
of all possible walks for $s=2$ when $p\in(0,1)$, plotted top-down to resemble a Galton board figure.
(The tree is truncated after $9$ steps, for more steps would only produce a  black band between 
the current cutoff and the finish line).
 Also shown (in red) is a computer-generated sample of $7$ random Riemann-$\zeta$ walks with $p=\frac13$ \&\ $s=2$.

\medskip

\hspace{275pt} \boxed{\textrm{{Fractal\ tree\ \&\ 7 walks)}}}\vspace{-15pt}

\includegraphics[scale=0.40]{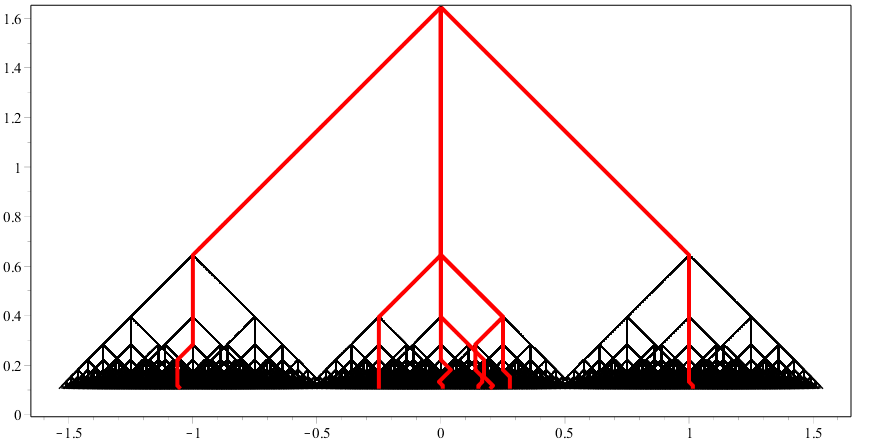}
\medskip

We also exhibit a histogram of the endpoints of $10^5$ walks with 1000 steps (Fig.~4).
\medskip

\hspace{270pt} \boxed{\textrm{{Histogram\ ($s=2$, $p=\tfrac13$)}}}\vspace{-15pt}

\includegraphics[scale=0.315]{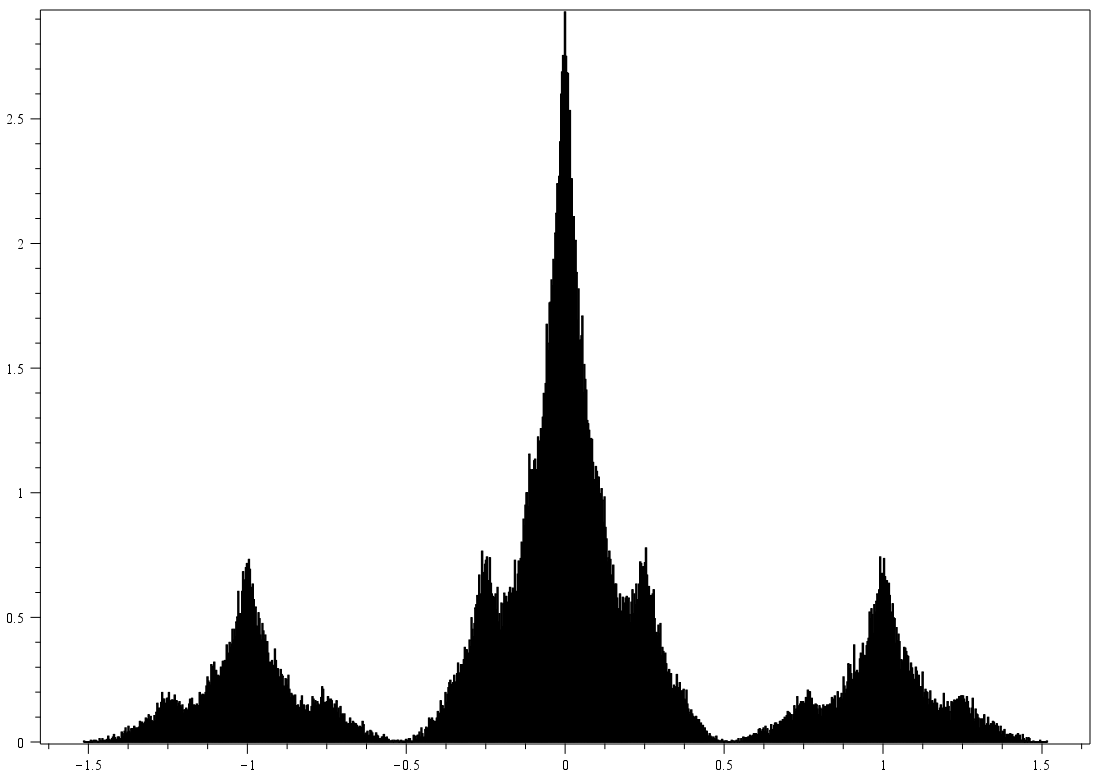}

 We next pick $s=1$ and two different choices of $p$, namely  $p=\frac13$ and $p=1$.
 For $s=1$ the random Riemann-$\zeta$ walks become so-called ``random Harmonic Series,'' which 
have been studied by Kac \cite{Kac}, Morrison \cite{Morrison}, and Schmuland \cite{Schmuland} 
in the special case that $p=1$.
 When $p\neq 1$ these harmonic random walks are interesting variations on their theme.
 We refrain, though, from trying to display the infinitely long harmonic random walk tree, for it is difficult to 
illustrate it faithfully.
 Yet the histograms of the endpoints of  $10^5$ harmonic walks with $10^3$ steps when $p=\frac13$ 
(Fig.~5) and $p=1$ (Fig.~6) are easily generated.
\bigskip

\hspace{270pt} \boxed{\textrm{{Histogram\ ($s=1$, $p=\tfrac13$)}}}\vspace{-15pt}

 \includegraphics[scale=0.3]{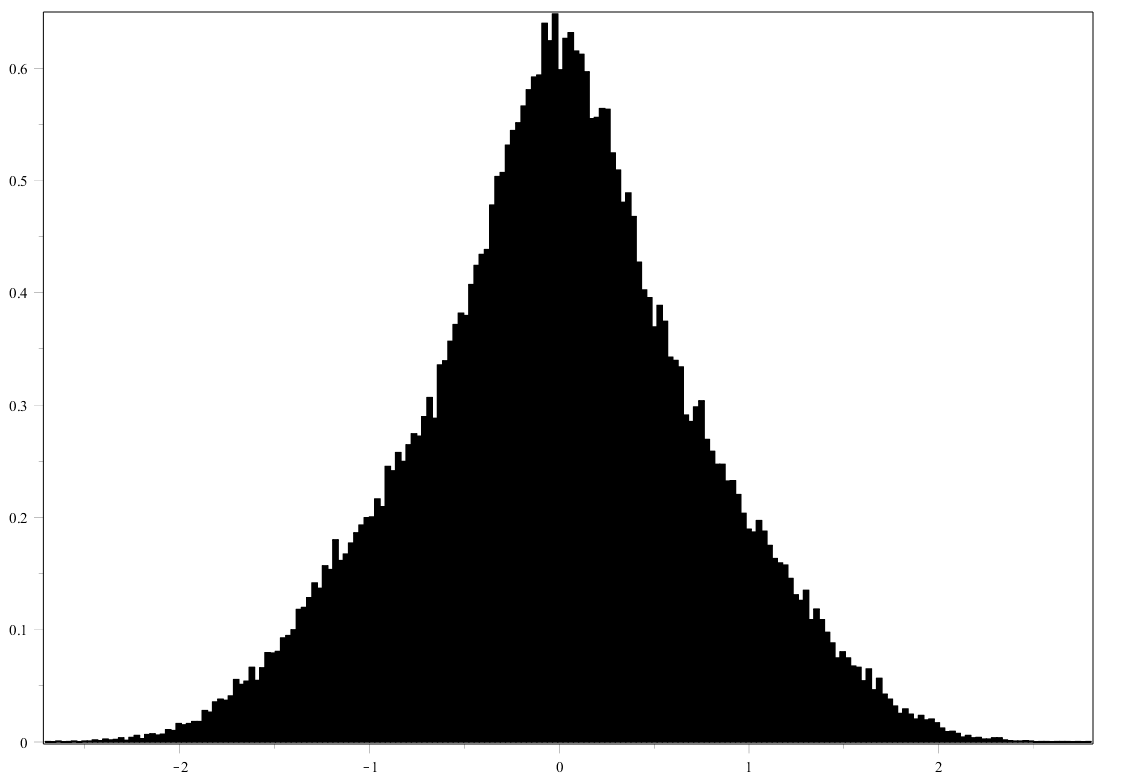}

\noindent

\hspace{270pt} \boxed{\textrm{{Histogram\ ($s=1$, $p=1$)}}}\vspace{-15pt}

\includegraphics[scale=0.3]{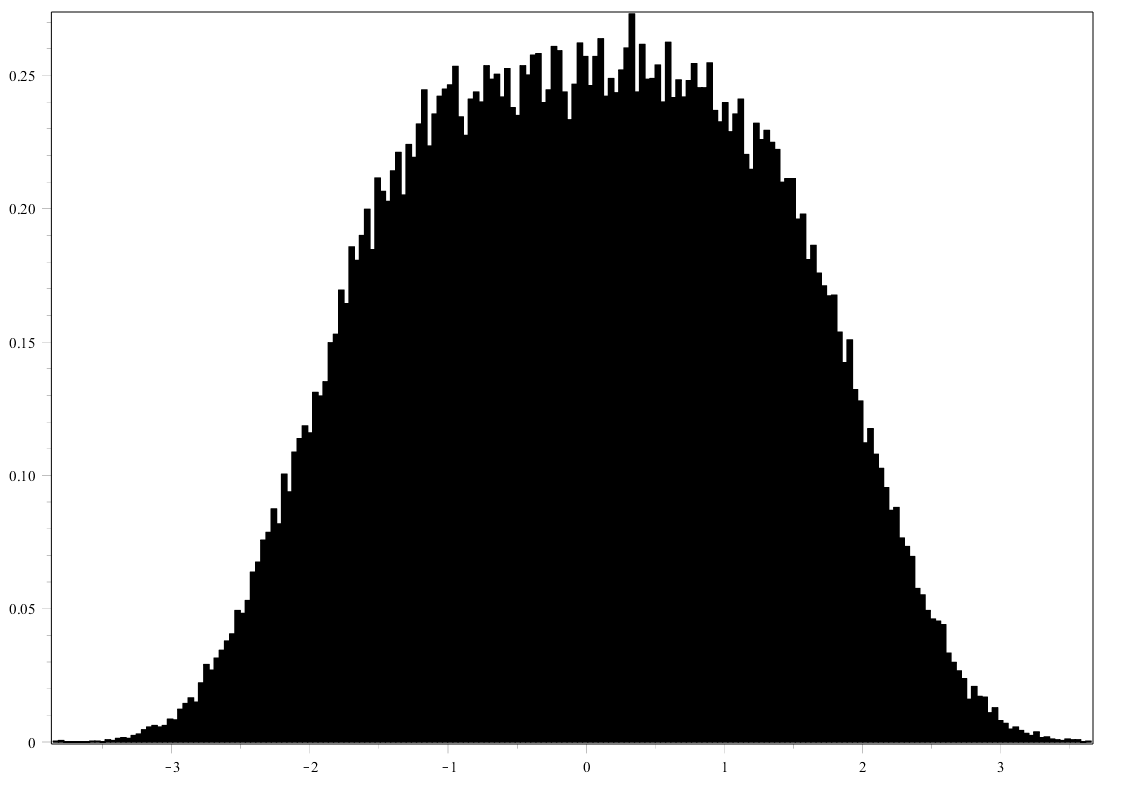}

\noindent

 Our Fig.~6 resembles the smooth theoretical PDF of the endpoints of the harmonic walk with $p=1$, displayed in 
Fig.~3 of \cite{Morrison} and Fig.~1 of  \cite{Schmuland}, quite closely; cf. the histogram based on $5,000$ walks 
with 100 steps displayed in Fig.~4 of \cite{Morrison}.
 When $p=1$ one is always on the move, so the histogram is quite broad.  
 Our Fig.~5 indicates that reducing $p$ (in this case to $p=1/3$) will lead to the build-up of a ``central peak.''
 The peak is even more pronounced in our Fig.~4 (where $p=1/3$ and $s=2$) which reveals a rich, conceivably 
self-similar structure with side peaks, and side peaks to the side peaks. 
 Our Fig.~4 also makes one wonder whether the peaks, if not fractal, could indicate that
 the first or second derivative of a theoretical PDF might blow up. 
 These questions  will be investigated in section \ref{THM}.

 But first, after having introduced random Riemann-$\zeta$ functions, at this point it is appropriate to 
pause and explain their relationship with the Riemann hypothesis.

\section{Typicality and the Riemann Hypothesis}\label{typical}

 Loosely speaking, a \emph{typical feature} of a random Riemann-$\zeta$ walk is a feature which ideally occurs 
``with probability 1'' (strong typicality), or at least ``in probability'' (weak typicality); see below. 
 A \emph{(strongly or weakly) perfectly typical random Riemann-$\zeta$ walk} is an empirical outcome of a 
random Riemann-$\zeta$ function evaluated at $s$
which \emph{exhibits all (either strongly or weakly) typical features}. 
 
 Since coin tosses are involved, for simplicity we look at the example of the set of all infinitely long sequences of 
fair coin tosses first. 

\subsection{Typicality for coin toss sequences}\label{typicalCOINS}

 We identify the events $H$ with $1$ and $T$ with $0$, and introduce the Bernoulli random variable $B\in\{0,1\}$,
with Prob$(B=1)=\frac12$. 
 Let $B_n$ be an identical and independent copy of $B$. 
 Then by the \emph{strong law of large numbers} (see \cite{probBOOK}) one has
\begin{equation}
\mbox{Prob}\left(\lim_{N\to\infty}\frac{1}{N}{\textstyle\sum\limits_{n=1}^N} B_n = \frac12\right) =1
\end{equation}
whereas the  \emph{weak law of large numbers} (see \cite{probBOOK}) says that for any $\epsilon>0$, 
\begin{equation}
\lim_{N\to\infty} \mbox{Prob}\left(\left|\frac{1}{N}{\textstyle\sum\limits_{n=1}^N} B_n - \frac12\right| > \epsilon\right) =0.
\end{equation}
 Let $b_n^{}\in\{0,1\}$ denote the outcome of the coin toss $B_n$. 
 Then based on either the strong, or the weak law of large numbers we say that 
``$\lim_{N\to\infty}\frac{1}{N}\sum_{n=1}^N b_n = \frac12$'' is a \emph{strongly, or weakly, typical feature} for such an 
empirical sequence of outcomes $\{b_n^{}\}_{n\in\Nset}^{}$. 
 Of course, not every empirical sequence $\{b_n^{}\}_{n\in\Nset}^{}$ does exhibit this typical feature; take, for instance,
$\{b_n^{}\}_{n\in\Nset}^{} =\{1,1,1,1,...\}$. 
 We therefore say that $\{1,1,1,1,...\}$ is \emph{an atypical empirical sequence} for the fair coin tossing process. 
 More generally, \emph{any} empirical sequence $\{b_n^{}\}_{n\in\Nset}^{}$ for which 
$\Big|\frac{1}{N}\sum_{n=1}^N b_n - \frac12\Big| > \epsilon$ occurs infinitely often is said to be an \emph{atypical
empirical sequence} for this coin tossing process.  

 Next, consider the sequence $\{b_n^{}\}_{n\in\Nset}^{} =\{1,0,1,0,1,0,...\}$. 
 Could this be a perfectly typical sequence? 
 Clearly, $\lim_{N\to\infty}\frac{1}{N}\sum_{n=1}^N b_n = \frac12$, 
but anyone who has ever flipped a coin a dozen times, again and again, knows that ``typically'' it doesn't happen to obtain
six consecutive 
1-0 pairs --- here we borrow the common sense usage of ``typicality;'' indeed, on average the alternating 
pattern of six consecutive 1-0 pairs occurs less than once in 4,000 repetitions of a dozen coin tosses, and the likelihood of
$k$ 1-0 pairs decreases to zero with $k$ increasing to infinity in a trial of length $2k$.

 Yet, in an infinite sequence of coin tosses, with probability 1 the pattern of six consecutive 1-0 pairs occurs infinitely often;
more generally, for \emph{any} $k\in\Nset$, with probability 1 a pattern with $k$ consecutive 1s, or a pattern with $k$ consecutive 
0s, as well as $k$ consecutive 1-0 pairs, all occur infinitely often.
 Thus \emph{recurrences of such $k$-patterns are strongly typical features} of this coin tossing process.

 Let's look at one more {strongly typical feature} --- a variation on this theme will turn out to be related to the Riemann hypothesis.
 Namely, since by either the weak or the strong law of large numbers we can informally say that when $N$ is large enough
then $\sum_{n=1}^N b_n \approx \frac12 N$, i.e. $\sum_{n=1}^N (2b_n-1) \approx 0$ in a perfectly typical empirical sequence, 
we next ask for the typical size of the fluctuations about this theoretical mean, i.e. how large can they be, typically?
 \emph{Khinchin's law of the iterated logarithm} states that for any $\epsilon > 0$, 
with probability 1 the event $|\sum_{n=1}^N (2B_n -1)| > (1-\epsilon) \sqrt{2N\ln\ln N}$ will occur infinitely often,
while the event 
$|\sum_{n=1}^N (2B_n -1)| > (1 +\epsilon) \sqrt{2N\ln\ln N}$ has  probability 0 of occurring infinitely often
in the sequence $\{B_n\}_{n\in\Nset}^{}$.
 Thus, 
\begin{equation}
\Big|\sum_{n=1}^N (2b_n -1)\Big| > (1-\epsilon) \sqrt{2N\ln\ln N}
\end{equation}
occurs for infinitely many $N$ in a perfectly typical empirical sequence $\{b_n^{}\}_{n\in\Nset}^{}$, and
\begin{equation}
\Big|\sum_{n=1}^N (2b_n -1)\Big| > (1+\epsilon) \sqrt{2N\ln\ln N}
\end{equation}
will happen at most finitely many times.
 
 Countlessly many more features occur with probability 1, many of them trivially 
(like Prob$(\sum_{n=1}^NB_n<N+\epsilon)=1$), but many others not, and some of them are deep.
 This makes it plain that it is impossible, or at least
extremely unlikely, that anyone will ever give an \emph{explicit} characterization of a \emph{perfectly typical empirical sequence} 
of coin tosses. 
 (It is even conceivable that no such sequence exists!)
 By contrast, once a particular feature has been proven to occur with probability 1 (the strong version), 
or in probability (the weak version), it is straightforward to ask whether a given empirical sequence exhibits this
particular \emph{aspect of typicality}. 

 We are now armed to address the connection of the  Riemann hypothesis with the notion of typicality of random Riemann-$\zeta$
functions.

\subsection{Typicality for random Riemann-$\zeta$ functions}\label{typicalRH}

 We begin by listing a few typical features of random Riemann-$\zeta$ walks. 
 Let $r_p^{}(n)\in\{-1,0,1\}$ denote the outcome of the random variable $R_p^{}(n)$, and for given $p\in(0,1]$ and $s>\frac12$ 
let $\omega_{p}^\zeta(s)$ denote the outcome for the random Riemann-$\zeta$ walk $\Omega_{p}^\zeta(s)$, i.e.
\begin{equation}\label{RZwOUT}
 \omega_{p}^\zeta(s) = \left\{{\textstyle\sum\limits_{n=1}^N} r_p^{}(n)\frac{1}{n^s}\right\}^{}_{N\in\Nset}.
\end{equation}
 Then the fair coin tossing process of the previous subsection now yields that
\begin{equation}\label{RZwTYPa}
\lim_{N\to\infty}\textstyle{\frac1N\sum\limits_{n=1}^N} r_p^{}(n) =0
\end{equation}
is a feature typically exhibited by an outcome $\omega_{p}^\zeta(s)$, independently of $p$ and $s$.
 Next, 
\begin{equation}\label{RZwTYPb}
\lim_{N\to\infty}\textstyle{\frac1N\sum\limits_{n=1}^N} |r_p^{}(n)| = p 
\end{equation}
is a $p$-dependent feature typically exhibited by an $\omega_{p}^\zeta(s)$, independently of $s$.
 Lastly, Rademacher's result mentioned above actually shows that typically
\begin{equation}\label{RZwTYPc}
\lim_{N\to\infty} \textstyle{\sum\limits_{n=1}^N} r_p^{}(n)\frac{1}{n^s} = \omega_{p}^\zeta(s) 
\end{equation}
exists on the real $\omega$ line whenever $s>\frac12$.
 All these are strongly typical features.

 We now inquire into the typicality of the following outcomes of random Riemann-$\zeta$ functions with $s>\frac12$:
Riemann's $\zeta$-function (\ref{zetaFCT}) itself, viz. $\zeta(s)=\sum_{n\in\Nset}1/n^s$ understood as a (not necessarily
convergent) sequence of its partial sums; its reciprocal 
\begin{equation}\label{reciZETAfct}
\frac{1}{\zeta(s)}=\textstyle{\sum\limits_{n\in\Nset}}\mu(n)\frac{1}{n^s},
\end{equation}
where $\mu(n) \in\{-1,0,1\}$ is the M\"obius function (see \cite{Edwards}); and also
\begin{equation}\label{ratioZETAfcts}
\frac{\zeta(2s)}{\zeta(s)}
=\textstyle{\sum\limits_{n\in\Nset}}\lambda(n)\frac{1}{n^s},
\end{equation}
where $\lambda(n) \in\{-1,1\}$ is Liouville's $\lambda$-function (see \cite{OEIS}).
 All are possible outcomes of a random Riemann-$\zeta$ walk with $s>\frac12$, any\footnote{$\zeta(s)$ and $\zeta(2s)/\zeta(s)$ are
possible outcomes also when $p=1$, while $1/\zeta(s)$ is not.} $p\in(0,1)$.
 In terms of the outcomes $r_p^{}(n)$ of the coin tossing process, Riemann's zeta function corresponds to 
 $r_p^{}(n)=1$ for all $n\in\Nset$, its reciprocal to  $r_p^{}(n)=\mu(n)$, and the ratio
${\zeta(2s)}/{\zeta(s)}$ to $r_p^{}(n)=\lambda(n)$. 
 Can any of these $\omega_{p}^\zeta(s)$ be perfectly typical outcomes, at least for some $p$ values?

 As to $\zeta(s)$ itself, it is clear that it must be atypical, since $r_p^{}(n)=1$ for all $n\in\Nset$
manifestly violates the $p$- and $s$-independent typicality feature (\ref{RZwTYPa}).
 Yet $\zeta(s)$ does not necessarily violate each and every aspect of typicality! 
 For instance, if $p=1$ then (\ref{RZwTYPb}) holds for $\zeta(s)$ (though not for any other $p\in(0,1)$).
 Moreover, while the sequence of its partial sums diverges to infinity when $s\in(\frac12,1]$ in violation of
the typicality feature (\ref{RZwTYPc}), this feature is verified by $\zeta(s)$ if $s>1$. 
 In any event, since $\zeta(s)$ is an \emph{extreme outcome}, it is intuitively clear that it will violate most
aspects of typicality --- in this sense, we say that $\zeta(s)$ is \emph{extremely atypical} for all $p\in(0,1]$.

 On to its reciprocal. 
 It is known that the \emph{Prime Number Theorem}\footnote{This is the statement that the number of primes
less than $x$ is asymptotically given by $\int_2^x\frac{d\xi}{\ln\xi}$, with relative error going
to zero as $x\to\infty$; see  \cite{Edwards}.}
is equivalent to the actual frequencies of the values $\mu(n)=1$ and $\mu(n)=-1$ being equal in the long run,
so ${1}/{\zeta(s)}$ exhibits the typicality feature (\ref{RZwTYPa}).
 It is also known that the actual frequency of values $\mu(n)\neq 0$ equals $1/\zeta(2)\; (=6/\pi^2)$ 
in the long run, so ${1}/{\zeta(s)}$ also exhibits the typicality feature (\ref{RZwTYPb}) if $p=1/\zeta(2)$ (though clearly not
for any other $p$ value).
 Furthermore, ${1}/{\zeta(s)}$ satisfies the typicality feature (\ref{RZwTYPc}) for all $s>\frac12$.
 Could ${1}/{\zeta(s)}$ perhaps be a {perfectly typical} random Riemann-$\zeta$ function for all $s>\frac12$
when $p=1/\zeta(2)$?
 Recall that this would mean that for each $s>\frac12$ the pertinent actual walk is a \emph{perfectly typical walk},
i.e. a walk which \emph{exhibits all features} of the theoretical random-walk law \emph{which occur with probability} 1 
(or at least \emph{in probability}).

 Similarly, the Prime Number Theorem is equivalent to the actual frequencies of the values $\lambda(n)=1$ and $\lambda(n)=-1$ 
being equal in the long run, so also the ratio ${\zeta(2s)}/{\zeta(s)}$ exhibits the typicality feature (\ref{RZwTYPa}).
 Furthermore, if (and only if) $p=1$ then ${\zeta(2s)}/{\zeta(s)}$ exhibits the typicality feature (\ref{RZwTYPb}).
 Lastly, ${\zeta(2s)}/{\zeta(s)}$ also exhibits the typicality feature (\ref{RZwTYPc}) for all $s>\frac12$.
 Could also ${\zeta(2s)}/{\zeta(s)}$ perhaps be a {perfectly typical} random Riemann-$\zeta$ function for all $s>\frac12$
when $p=1$?

 A moment of reflection reveals that this would be \emph{truly paradoxical}: if 
${1}/{\zeta(s)}$ and / or ${\zeta(2s)}/{\zeta(s)}$ are {perfectly typical} random Riemann-$\zeta$ functions for the
mentioned $p$-values, then one can learn a lot about them by studying what is {typical} for random walks with those $p$-values,
without ever looking at ${1}/{\zeta(s)}$ or ${\zeta(2s)}/{\zeta(s)}$.
 Of course, if one learns something about ${1}/{\zeta(s)}$ and / or ${\zeta(2s)}/{\zeta(s)}$, then one also 
learns something about $\zeta(s)$ --- but how can one learn something about an extremely atypical random Riemann-$\zeta$ function by
studying what is typical for such random walks? 
 The obvious way out of this dilemma is to conclude:\hfill

\centerline{\emph{Neither ${1}/{\zeta(s)}$ nor ${\zeta(2s)}/{\zeta(s)}$ can be perfectly typical random Riemann-$\zeta$ functions!}}
\newpage

 The upshot is that both ${1}/{\zeta(s)}$ and ${\zeta(2s)}/{\zeta(s)}$ \emph{must feature} {some} \emph{atypical empirical statistics}, 
encoded in the sequences $\{\mu(n)\}^{}_{n\in\Nset}$ and $\{\lambda(n)\}^{}_{n\in\Nset}$.
 Obviously these atypical features must be inherited from the correlations in the distribution of prime numbers;
recall that the coin tossing process, by contrast, is correlation-free.
 Since the Riemann hypothesis about the location of the non-real zeros of $\zeta(s)$ is equivalent to some detailed knowledge about
the distribution of and correlations amongst prime numbers, it may well be that some particular atypical empirical feature of 
${1}/{\zeta(s)}$ and ${\zeta(2s)}/{\zeta(s)}$ will be equivalent to the Riemann hypothesis. 
 Which kind of feature, if any, remains anybody's best guess  --- to the best of our knowledge.

 Surprisingly, and indeed intriguingly, it is known though that \emph{a certain typical feature},
if indeed exhibited by the $1/\zeta(s)$ walk for $p=1/\zeta(2)$, beyond the agreement of empirical and theoretical frequencies,
\emph{is equivalent to the Riemann hypothesis}!
 We are grateful to Alex Kontorovich for having pointed this out to us.

 Namely, let us extend the definition of the random Riemann-$\zeta$ walk $1/\zeta(s)$ to $s=0$, 
\emph{not} by analytic continuation, but in terms of the sequence of its partial sums:
\begin{equation}\label{reciRZwATnull}
\frac{1}{\zeta(0)} := \left\{\textstyle\sum\limits_{n=1}^N \mu(n)\right\}^{}_{N\in\Nset}.
\end{equation}
 Note that for $s\leq \frac12$ the $1/\zeta(s)$ random walk may well wander off to infinity, but the rate at which this happens
is crucial (recall Khinchin's law of the iterated logarithm which we mentioned in subsection \ref{typicalCOINS}).
 As explained in \cite{Edwards}, chpt.12.1, Littlewood proved the equivalence: 
\begin{equation}\label{reciRZwATnullGROWTH}
\forall\epsilon>0:\ 
\lim_{N\to\infty}N^{-\frac12-\epsilon}\Big|\sum_{n=1}^N\mu(n)\Big| =0\quad \leftrightarrow\quad
 \mbox{The\ Riemann\ hypothesis\ is\ true}.
\end{equation}
 And as explained in \cite{Edwards}, chpt.12.3, Denjoy noted that if one assumes that the $\pm1$ values of $\mu(n)$ are 
distributed as if they were generated by fair and independent coin flips, then the central limit theorem implies that 
$\lim_{N\to\infty}N^{-\frac12-\epsilon}|\sum_{n=1}^N\mu(n)| =0$ holds with probability 1.
 Of course, $\mu(n)=0$ is still determined by its formula, but the empirical frequency of $\mu(n)=0$ occurrences
is $1-6/\pi^2$ in the long run, and by adopting Denjoy's reasoning one can show that for $p= 6/\pi^2$ one has that 
\begin{equation}\label{randomDenjoy}
\forall\epsilon>0:\ \mbox{Prob}\left(\lim_{N\to\infty}N^{-\frac12-\epsilon}\Big|\sum_{n=1}^N R^{}_{6/\pi^2}(n)\Big| =0\right) =1. 
\end{equation}
 Thus l.h.s.(\ref{reciRZwATnullGROWTH})  would be a typical feature exhibited by the $1/\zeta(0)$ walk at $p= 6/\pi^2$.
 \medskip


\section{The Characteristic Function of $\Omega^\zeta_p(s)$}\label{PROB}
  We now show that the infinite trigonometric products $\Clg(t)$ given in (\ref{genCLOITREprod}) 
are characteristic functions of the $\Omega_{p}^\zeta(s)$, i.e. $\Clg(t) = \mbox{Exp}\big(\exp\big(it\Omega_{p}^\zeta(s)\big)\big)$,
where ``Exp'' is \emph{expected value} (not to be confused with the exponential function $\exp$).
 Since $\Omega_{p}^\zeta(s)$ is an infinite sum of independent random variables $R_p^{}(n)/n^s$ (see (\ref{RANDOMzetaFCT})), 
$\exp\!\big(it\Omega_{p}^\zeta(s)\big)$ is an infinite product of independent random variables $\exp\!\big(itR_p^{}(n)/n^s\big)$, and
by a well-known theorem in probability theory, expected values of products of independent random variables are products of the
their individual expected values. 
 And so we have (temporarily ignoring issues of convergence)
\begin{eqnarray}\label{charFUNC}
\mbox{Exp}\Big(\exp\big(it\Omega_{p}^{\zeta}(s)\big)\Big)
\!\!&=&\!\!  \mbox{Exp}\Big(\prod_{n\in\Nset}\exp\big(itR_{p}^{}(n)\tfrac{1}{n^s}\big)\Big)
=  \prod_{n\in\Nset}\mbox{Exp}\Big(\exp\big(itR_{p}^{}(n)\tfrac{1}{n^s}\big)\Big)\cr
\!\!&=&\!\! \prod_{n\in\Nset} \Big(\tfrac{1}{2}p e^{-i\hfrac{t}{n^s}}
+ (1-p) + \tfrac{1}{2}p e^{i\hfrac{t}{n^s}}\Big) \cr
\!\!&=&\!\! \prod_{n\in\Nset} \Big(1-p + p \cos\big(\tfrac{t}{n^s}\big)\Big)
 \equiv \Clg(t),\vspace{-15pt}
 \end{eqnarray}
where we have used Euler's formula to rewrite $\tfrac12\big(e^{i\hfrac{t}{n^s}} +e^{-i\hfrac{t}{n^s}}\big)=\cos\big(\hfrac{t}{n^s}\big)$.

 That was straightforward.
 Next we explain the relationship between the characteristic functions $\Clg(t)$ of $\Omega^\zeta_p(s)$ 
and the probability distribution $\varrho^{\zeta}_{p;s}(d\omega)$ of the endpoints of these 
random Riemann-$\zeta$ walks on the $\omega$ line.
 Formally this is accomplished by realizing that $\Clg(t)$ is the inverse Fourier transform of $\varrho^{\zeta}_{p;s}(d\omega)$, viz.
\begin{eqnarray}\label{charFUNCisFOURIERinverseOFdistribution}
\mbox{Exp}\Big(\exp\big(it\Omega_{p}^{\zeta}(s)\big)\Big)
= \int_{\Rset} e^{it\omega}\varrho^{\zeta}_{p;s}(d\omega).
 \end{eqnarray}
 Therefore we obtain $\varrho^{\zeta}_{p;s}(d\omega)$ by taking the Fourier transform of $\Clg(t)$. 
  As recalled in \cite{Morrison}, the Fourier transform of a product equals the \emph{convolution product} (``$*$'', 
see below) of the Fourier transforms of its factors, and so we find
\begin{eqnarray}\label{RZdistribution}
\varrho^{\zeta}_{p;s}(d\omega)
\!\!&=&\!\! \Big(*\prod_{n\in\Nset} \frac{1}{2\pi}\int_\Rset e^{-i\omega t} 
\Big[\tfrac{1}{2}p e^{-i\hfrac{t}{n^s}}
+ (1-p) + \tfrac{1}{2}p e^{i\hfrac{t}{n^s}}\Big]dt\Big) (d\omega)\cr
\!\!&=&\!\! \Big(*\prod_{n\in\Nset} \Big[\tfrac{1}{2}p\delta_{-\frac{1}{n^s}}^{} 
+ (1-p)\delta_{0}^{} 
+ \tfrac{1}{2}p \delta_{\frac{1}{n^s}}^{}\Big]\Big)(d\omega);
\end{eqnarray}
here, ``$*\prod$'' denotes repeated convolution (cf. \cite{Morrison}), and
$\delta_{\omega_k}^{}$ is a Dirac measure.\footnote{If
    $I\subset\Rset$ is any closed interval, then $\int_I \delta_{\omega_k}^{} (d\omega)= 1$ if 
    $\omega_k^{}\in I$ and $\int_I \delta_{\omega_k^{}}^{} (d\omega)= 0$ if $\omega_k^{}\not\in I$.}

 This distribution looks intimidating, but it only conveys what we know already!
 Namely, formally (\ref{RZdistribution}) is the limit $N\to\infty$ of the $N$-fold partial convolution products\footnote{We temporarily 
   suppress the suffix ``$\zeta$'' so as not to overload the notation.} 
\begin{eqnarray}\label{NstepsDELTAconcise}
\varrho^{(N)}_{p;s}(d\omega)
:= \Big(*\prod_{n=1}^N \Big[\tfrac{1}{2}p\delta_{-\frac{1}{n^s}}^{} 
+ (1-p)\delta_{0}^{} 
+ \tfrac{1}{2}p \delta_{\frac{1}{n^s}}^{}\Big]\Big)(d\omega).
\end{eqnarray}
 Now recall that the {convolution} product, which for two integrable functions $f$ and $g$ 
is defined by $(f*g)(\omega) = \int f(\omega')g(\omega-\omega')d\omega'$, extends
to delta measures where it acts as follows: $\delta_a^{}*\delta_b^{} = \delta_{a+b}^{}$ (see \cite{Morrison}).
 Therefore, by multiplying out the convolution product at r.h.s.(\ref{NstepsDELTAconcise}), 
using the distributivity of ``$*$'' one finds that $\varrho^{(N)}_{p;s}(d\omega)$ is a
weighted sum of point measures at the possible outcomes
\vspace{-20pt}

\begin{equation}\label{RZwNout}
 \omega_{p}^{(N)}(s)
: = \sum_{n=1}^N r_p^{}(n)\frac{1}{n^s} \in {\cal{L}}_{p}^{(N)}(s),\quad r_p^{}(n) \in \{-1,0,1\}, \vspace{-5pt}
\end{equation}
of the random walk truncated after $N$ steps, 
\vspace{-20pt}

\begin{equation}\label{RZwN}
 \Omega_{p}^{(N)}(s)
 := \sum_{n=1}^N R_p^{}(n)\frac{1}{n^s}.
\end{equation}

\vspace{-15pt}\noindent
 The set of locations ${\cal{L}}_{p}^{(N)}(s)\subset\Rset$ is finite, and generically\footnote{It 
  may in principle happen for certain discrete values of $s$ (but not of $p$) that different 
  $N$-step paths lead to the same outcome $\omega_{p}^{(N)}(s)$.
  However, since $s>\frac12$ is a continuous parameter, this degenerate situation is not generic.
  Note though that it may well happen that we humans ``inadvertendly'' pick precisely those non-generic cases,
  for instance if degeneracy occurs when $s\in\Nset$!}
consists of $3^N$ distinct real points if $p\in(0,1)$, and of $2^N$ distinct real points if $p=1$. 
 Thus, $\varrho^{(N)}_{p;s}(d\omega)$ becomes
\begin{equation}\label{NstepsDELTA}
\varrho^{(N)}_{p;s}(d\omega) 
= \sum_{{\cal{L}}_{p}^{(N)}(s)} \mbox{P}\big(\omega_p^{(N)}(s)\big)\delta_{\omega_p^{(N)}(s)}^{} (d\omega);
\end{equation}

\vspace{-10pt}\noindent
the sum runs over all ${\omega_p^{(N)}(s)\in {\cal{L}}_{p}^{(N)}(s)}$, and 
$\mbox{P}\big(\omega_p^{(N)}(s)\big):=\mbox{Prob}\big( \Omega_{p}^{(N)}(s)=\omega_p^{(N)}(s)\big)$.
 These probabilities $\mbox{P}\big(\omega_p^{(N)}(s)\big)$ are readily computed from the tree diagram in Fig.~3,
or by inspecting (\ref{RZwNout}):
 if in order to reach $\omega_p^{(N)}(s)$ you need to move $m\leq N$ times (whether left or right has equal probability),
then $\mbox{P}\big(\omega_p^{(N)}(s)\big) = (p/2)^m (1-p)^{N-m}$, independently of $s$.
 Note that there are $2^m \genfrac{(}{)}{0pt}{}{N}{m}$ possible outcomes for an $N$-step walk with $m\leq N$ 
moves, and indeed $\sum_{m=0}^N 2^m \genfrac{(}{)}{0pt}{}{N}{m}(p/2)^m (1-p)^{N-m} = (1-p+p)^N=1$.

 Let's look at two examples.
 After $1$ step with $p\in(0,1)$ there are 3 possible positions, and the distribution (\ref{NstepsDELTAconcise}) with $N=1$
reads
\begin{equation}\label{ONEstepsDELTA}
\varrho^{(1)}_{p;s}(d\omega)
= \Big(\tfrac{1}{2}p \delta_{-1}^{}+ (1-p) \delta_0^{} + \tfrac{1}{2}p \delta_{1}^{}\Big)(d\omega) .
\end{equation}
 After $2$ steps with $p\in(0,1)$ we have 9 possible positions, and (\ref{NstepsDELTAconcise})  with $N=2$ reads
\begin{eqnarray}\label{TWOstepsDELTA}
\varrho^{(2)}_{p;s}(d\omega)
=& \Big(\Big[\tfrac12 p \delta_{-1}^{} + (1-p)\delta_0^{} +\tfrac12 p\delta_{1}^{}\Big]*\Big[\tfrac{1}{2}p\delta_{-\frac{1}{2^s}}^{} 
+ (1-p)\delta_{0}^{} + \tfrac{1}{2}p \delta_{\frac{1}{2^s}}^{}\Big]\Big)(d\omega)\cr
= &\Big( \tfrac{1}{4}p^2 \delta_{-1-\frac{1}{2^s}}^{} 
+ \tfrac{1}{2}p(1-p)\delta_{-1+0}^{} 
+ \tfrac{1}{4}p^2 \delta_{-1+\frac{1}{2^s}}^{} \Big)(d\omega) + \cr
&
 \Big(\tfrac{1}{2}p(1-p) \delta_{0-\frac{1}{2^s}}^{} 
+(1-p)^2 \delta_{0+0}^{}
+ \tfrac{1}{2}p(1-p)  \delta_{0+\frac{1}{2^s}}^{}\Big)(d\omega)
+\cr
&
\Big( \tfrac{1}{4}p^2\delta_{1-\frac{1}{2^s}}^{}
 +\tfrac{1}{2}p(1-p)\delta_{1+0}^{} 
 + \tfrac{1}{4}p^2 \delta_{1+\frac{1}{2^s}}^{} \Big)(d\omega),
\end{eqnarray}
which is precisely (\ref{NstepsDELTA}) with $N=2$;
we have facilitated the comparison by writing all two-step walks which lead to the locations of the point masses explicitly, including
the ``non-moves.''
 Similarly one can compute the $N$-th partial convolution product, although this soon gets cumbersome and does not illuminate
the process any further.

 The theory of convergence of probability measures (e.g. ref.[1] in \cite{Schmuland}) shows that the
sequence of partial products (\ref{NstepsDELTAconcise}) does converge to a probability measure (\ref{RZdistribution}) 
if $s>\frac12$.
 Unfortunately, the expression (\ref{RZdistribution}) does not readily give up its secrets.

 In particular, each measure (\ref{NstepsDELTAconcise}) is singular with respect to (w.r.t.) Lebesgue measure $d\omega$, 
so could it be that the $N\to\infty$ limit (\ref{RZdistribution}) is singular as well --- e.g., supported by a 
fractal?
 And if not, when $\varrho^{\zeta}_{p;s}(d\omega)$ is absolutely continuous w.r.t. $d\omega$, 
is its PDF perhaps not differentiable at its peaks, as hinted at by Fig.~4?
 
 The answers to these questions will be extracted from $\Clg(t)$ in the next section.
\vspace{-27pt}

\section{The Main Theorem}\label{THM}
\vspace{-10pt}

 In this section we use elementary calculus techniques to prove the following result:
%
%
%
\begin{theo}\label{THMgenCLOITREtrend}
 Let $p\in(0,1]\ \&\ s> \frac12$. 
 Then
\vspace{-2pt}
\begin{equation}\label{genClTHM}
\forall\,t\in\Rset:\quad 
{\mathrm{Cl}}_{p;s}^{}(t) = 
\exp\left({- C_{p;s}^{} \,|t|^{1/s}}\right) F_{p;s}^{}(|t|),
\end{equation}
where $|F_{p;s}^{}(|t|)| \leq \exp(K_{p;s}^{} |t|^{1/(s+1)})$ for some constant $K_{p;s}^{}>0$, and where
\begin{equation}\label{genCLOITREcoeffC}
C_{p;s}^{} 
:= -\frac1s\int_0^\infty\ln|{1-p+p\cos\xi}|\frac{1}{\xi^{1+1/s}}{\rm{d}}\xi.
\end{equation}

\vspace{-5pt}\noindent
Moreover, when $p\in(0,\frac12)\ \&\ s> \frac12$ the stronger bound $|\ln F_{p;s}^{}(|t|)| \leq K_{p;s}^{} |t|^{1/(s+1)}$ 
holds; furthermore, we then have
$C_{p;s}^{}  = A_{s}^{} B_{p;s}^{}$, with
\vspace{-5pt}
\begin{equation}\label{genCLOITREcoeffA}
A_{s}^{} := \int_0^\infty {\sin\xi}\frac{1}{\xi^{1/s}}{\rm{d}}\xi = \textstyle\Gamma\left(1-\frac1s\right)\cos\left(\frac{\pi}{2s}\right)
\end{equation}
(where it is understood that $A_1^{}=\lim_{s\to1}\Gamma\bigl(1-\frac1s\bigr)\cos\bigl(\frac{\pi}{2s}\bigr)\ [=\frac{\pi}{2}]$), and
\begin{equation}\label{genCLOITREcoeffB}
B_{p;s}^{} 
:= 
\sum_{n=0}^\infty(-1)^n\left(\frac{p}{1-p}\right)^{n+1}\frac{1}{2^n}\sum_{k=0}^{\ceil{\frac{n-1}{2}}}
\begin{pmatrix} n\cr k \end{pmatrix} \frac{({1+n-2k})^{1/s}\hspace{-5pt}}{{1+n-k}\;}\qquad.
\end{equation}
\end{theo}
\newpage

\begin{rema}
 Recalling that $\ln|z| = \Re{e} \ln z$ for $z\in\Cset$, we conclude from (\ref{genCLOITREcoeffC}) that
\begin{equation}\label{genCLOITREcoeffCrealPART}
C_{p;s}^{} 
= -\frac1s\Re{e}\int_0^\infty\ln({1-p+p\cos\xi})\frac{1}{\xi^{1+1/s}}{\rm{d}}\xi,\quad p\in(0,1]\ \&\  s>\frac12,
\end{equation}
where the integral at r.h.s.(\ref{genCLOITREcoeffCrealPART}) is understood as
analytic continuation from $p\in(0,\frac12)$ (when $\ln(1-p+p\cos\xi)\in\Rset$) to $p\in(\frac12,1]$
(when $\ln(1-p+p\cos\xi)\in\Cset$).
  Larry Glasser and Norm Frankel (personal communications, Dec. 2016) have informed us that
this analytically continued $C_{p;s}^{}$, denoted $\widetilde{C}_{p;s}^{}$,
has been calculated in \cite{Glasser}             
to
\begin{equation}\label{CpsLARRY}
\widetilde{C}_{p;s}= -2\textstyle\Gamma\left(1-\frac1s\right)\cos\left(\frac{\pi}{2s}\right)
{\mathrm{Li}}_{1-\frac1s}\!\left(\sqrt{q^2-1}-q\right), \quad q= (1-p)/p;
\end{equation}
here $\mathrm{Li}_{a}\!(z)$ is a polylogarithm. 
 They also remarked that for $p=\frac12$ and $s=1$ one has
\begin{equation}\label{genCLOITREcoeffCpvSPECIAL}
C_{\frac12;1}^{} = \mbox{\sc{pv}}\!\int_0^\infty\!\frac{\sin\xi}{1+\cos\xi}\frac{1}{\xi}{\rm{d}}\xi
= \mbox{\sc{pv}}\!\int_0^\infty\frac{\tan\xi}{\xi}{\rm{d}}\xi =C_{1;1}^{},
\end{equation}
with $C_{1;1}^{}=\frac\pi2$; see \cite{BG} (here ``{\sc{pv}}'' means \emph{principal value}).
 So presumably 
\begin{equation}\label{genCLOITREcoeffCpvGENERAL}
C_{p;s}^{}\! =\! \mbox{\sc{pv}}\!\int_0^\infty\!\!\!\frac{p\sin\xi}{1-p+p\cos\xi}\frac{1}{\xi^{1/s}}{\rm{d}}\xi
\end{equation}
for $p\in[\frac12,1]$ and $s>\frac12$.
 Note that $p\mapsto C_{p;s}$ has a derivative singularity at $p=\frac12$.
\end{rema}

\begin{rema}
 Identity (\ref{genCLOITREcoeffCpvSPECIAL}) is a special case of the identity $C_{1;s}^{} = 2^{-1+1/s}C_{1/2;s}^{}$;
$s>\frac12$.
 Indeed, a half-angle identity and an obvious substitution of variables yields

\vspace{-15pt}
\begin{eqnarray}\label{Chalfs}
\int_0^\infty\!\!\ln(\tfrac12+\tfrac12\cos\xi)\frac{1}{\xi^{1+1/s}}{\rm{d}}\xi 
=\!\! 
\int_0^\infty\!\!\ln(\cos^2\tfrac{\xi}{2})\frac{1}{\xi^{1+1/s}}{\rm{d}}\xi
= 2^{1-\frac1s}\!\!\int_0^\infty\!\!\ln|\cos {\xi}|\frac{1}{\xi^{1+1/s}}{\rm{d}}\xi.
\end{eqnarray}

\vspace{-5pt}\noindent
But l.h.s.(\ref{Chalfs})=$-sC_{\hfrac12;s}^{}$ while r.h.s.(\ref{Chalfs})=$-s2^{1-\frac1s}C_{1;s}^{}$.
\end{rema}

\begin{rema}
 It is also easy to show that \emph{$\mbox{Cl}_{1/2;s}^{}(t) =  \mbox{Cl}_{1;s}^2(t/2)$}.
\hfill End of Remarks.
\end{rema}

 Cloitre's surmise follows from Theorem \ref{THMgenCLOITREtrend}.
 Indeed, we have the stronger result:
\begin{coro}
 For $p=1/3$ and $s=2$, Theorem \ref{THMgenCLOITREtrend}  reduces to
\begin{equation}\label{CLOITREasymp}
\forall\,t\in\Rset:\quad 
{\mathrm{P}}_{\!\mbox{\tiny\sc{Cl}}} (t)= e^{- C \,\sqrt{|t|} +\varepsilon(|t|)},
\end{equation}
\vspace{-15pt}

\noindent
with correction term bounded as $|\varepsilon(|t|)|\leq K |t|^{1/3}$ for some $K>0$, and with
\vspace{-5pt}
\begin{equation}\label{CLOITREcoeffC}
C= \int_\Rset\frac{\sin\xi^2}{2+\cos \xi^2}{\rm{d}}\xi = 
\sqrt{\frac\pi2}\sum_{n=0}^\infty(-1)^n\frac{1}{2^{2n+1}}\sum_{k=0}^{\ceil{\frac{n-1}{2}}}
\begin{pmatrix} n\cr k \end{pmatrix} \frac{\sqrt{1+n -2k}}{1+n-k};
\end{equation}
numerically, 
 $C = 0.319905585... \sqrt{\pi}\approx 0.32\sqrt{\pi}$. 
\end{coro}
\begin{rema}
 By a simple change of variables, and the fact that $\xi^2$ is even, we have
\begin{equation}
C = \int_{\Rset}\frac{\sin\xi^2}{2+\cos\xi^2}{\rm{d}}\xi = 
2 \int_0^\infty\frac{\sin\xi^2}{2+\cos\xi^2}{\rm{d}}\xi =  
 \int_0^\infty\frac{\sin\xi}{2+\cos\xi}\frac{1}{\xi^{1/2}}{\rm{d}}\xi  = C_{\frac13;2}^{} 
\end{equation}
\end{rema}

 Theorem \ref{THMgenCLOITREtrend} implies that $\Omega_{p}^\zeta(s)$ is a very regular random variable. 
 Namely, $\Clg(t)\in C^\infty$, and by (\ref{genClTHM}) the integral of $|t|^m\Clg(t)$ exists for any $m\in\{0,1,2,...\}$, 
so $\Clg(t)$ is a Schwartz function, and so the Fourier transform of $\Clg(t)$ is a Schwartz function;
see chpt. IX in \cite{ReedSimonI}.
 Also, $\Clg(0)=1$, so its Fourier transform has integral 1.
 We already know that its Fourier transform is positive.
 Thus we have
\begin{coro}\label{regularity}
 The random variable ${\Omega_{p}^\zeta(s)}^{}$ defined by its characteristic function (\ref{genCLOITREprod}) with 
$p\in(0,1]$ and $s>\frac12$ converges with probability 1. 
 It is continuous, having a probability density
$f_{\Omega_{p}^\zeta(s)}^{}(\omega)$ which is a (generally not real analytic) Schwartz function.
\end{coro}

 Corollary \ref{regularity} settles our questions concerning the distribution of $\Omega_{1/3}^{\zeta}(2)$.
 Despite the seemingly self-similar structure hinted at in Fig.~4, $\Omega_{1/3}^{\zeta}(2)$ cannot be supported on a fractal
subset of the $\omega$ line.
 Also, despite the appearance of singular peaks hinted at in Fig.~4, the PDF of $\Omega_{1/3}^{\zeta}(2)$ is infinitely often
continuously differentiable. 

 It remains to prove Theorem~\ref{THMgenCLOITREtrend}.
 To offer some guidance we explain our \textit{strategy}.

 \textbf{First of all}, we prove the stronger part of Theorem~\ref{THMgenCLOITREtrend}
concerning $p\in(0,\frac12)$.
 We then have $\Clg(t)>0$, so we can take its logarithm and obtain the infinite series
\vspace{-5pt}
\begin{equation}\label{logCloitrePRODgeneralized}
\ln \Clg(t)
= \textstyle{\sum\limits_{n\in\Nset}}
\ln \left[1- p +p\cos\left(\textstyle n^{-s}_{_{}}t\right)\right],\ t\in\Rset,
\end{equation}
with $\textstyle p\in(0,\frac12)\ \&\ s> \frac12$.
 We then follow the proof of Theorem~1 of \cite{KieJSPforBenJerryMichael} which establishes that if $s>1$, then 
for all $t\in\Rset$ one has
$\sum_{n\in\Nset}\sin\left(\textstyle n^{-s}_{_{}}t\right) = \alpha_s^{}\sign(t)|t|^{1/s} + \varepsilon(|t|)$, with
$\alpha_s^{} = \Gamma\bigl(1-\frac1s\bigr)\sin\bigl(\frac{\pi}{2s}\bigr)$ and $|\varepsilon(|t|)| \leq K_s^{} |t|^{1/(s+1)}$ 
for some $K_s^{}>0$.
 
 Note that by the reflection symmetry about $t=0$ of $\Clg(t)$ it suffices to consider $t>0$.
 Yet one needs to distinguish $0\leq t\leq t_s^{}$ and $t\geq t_s^{}$ for some $t_s^{}>0$.

 The \emph{near side} $0\leq t\leq t_s^{}$ will be estimated with the help of a Maclaurin expansion and turn out to be
subdominant.

 The \emph{far side} $t\geq t_s^{}$ will be handled by splitting the infinite series into two parts,
\begin{equation}
\textstyle{\sum\limits_{n\in\Nset}}(\cdots_n) 
=\label{eq:SdefSsplitSplusONEroot}
\textstyle{\sum\limits_{n=1}^{N_s(t/\tau)}} (\cdots_n) + \textstyle{\sum\limits_{n=N_s(t/\tau)+1}^{\infty} }(\cdots_n),
\end{equation}
where $(\cdots_n) = \ln \left[1- p +p\cos\left(\textstyle n^{-s}_{_{}}t\right)\right]$, and where
$N_s(t/\tau): = \lceil{(st/\tau)^{1/(s+1)}}\rceil$, with $\tau <\pi/2$.
 The first (finite) sum in (\ref{eq:SdefSsplitSplusONEroot}) will be shown to yield only a subdominant error bound. 
 The second (infinite) sum in (\ref{eq:SdefSsplitSplusONEroot}) will be interpreted as a Riemann sum approximation 
to an integral over the real line, the trend function, plus a subdominant error bound.
 We now outline this argument.

	Since $\tau < \pi/2$, when $t$ gets large
any two consecutive arguments $t/n^s$ and $t/(n+1)^s$ of the cosine functions will come to lie within one quarter period of cosine
whenever $n > \lceil{(st/\tau)^{1/(s+1)}}\rceil$.
	Moreover, with increasing $n$, for fixed $t/\tau$, the consecutive arguments $t/n^s$ and $t/(n+1)^s$ 
will be more and more closely spaced.
	And so, when $\tau$ is sufficiently small, with increasing $t$ the part of the sum of $\ln\Clg(t)$ with 
 $n > N_s(t/\tau)$ becomes an increasingly better Riemann sum approximation to 
$$
\int_{ N_s(t/\tau)+1}^\infty \ln \left[1- p +p\cos\left(\textstyle \nu^{-s}_{_{}}t\right)\right] {\rm{d}}\nu.
$$
	More precisely, using the variable substitution $\nu^{-s}t=\xi$, we have (informally)
\begin{equation}
\hspace{-25pt}	
{\textstyle\sum\limits_{n= N_s(t/\tau)+1}^{\infty}}\hspace{-10pt}\ln\left[1- p +p\cos\left(\textstyle n^{-s}_{_{}}t\right)\right]
\approx \label{eq:SdefSsplitSplusONErootINTapprox}
	t^{1/s}	{\textstyle\frac{1}{s}}\displaystyle
	\int_0^{t/( N_s(t/\tau)+1)^s} \hspace{-30pt}\ln\left[1- p +p\cos \xi \right]\frac{1}{\xi^{1+1/s}}  {\rm{d}}\xi.
\end{equation}
	Since $s>1/2$, the upper limit of integration at r.h.s.(\ref{eq:SdefSsplitSplusONErootINTapprox})
goes~to~$\infty$ like $K t^{1/(s+1)}$ when $t\to\infty$. 
 The limiting integral is an improper Riemann integral which converges absolutely for all $s>1/2$, yielding
\begin{equation}
	{\textstyle\frac1s}\int_0^\infty \ln\left[1- p +p\cos \xi \right]\frac{1}{\xi^{1+1/s}}  {\rm{d}}\xi
\equiv \label{eq:alphaSint}
	 - C_{p;s}^{}.
\end{equation}
 This heuristic argument will be made rigorous by supplying the subdominant error bounds, using only senior level undergraduate
mathematics.

 The integral (\ref{eq:alphaSint}) will be evaluated with the help of a rapidly converging geometric
series expansion and a recursion which involves the Catalan numbers.

 \textbf{Secondly}, we consider the regime $p\in[\tfrac12,1]$. 
 In this case $\Clg(t)$ has zeros at
\begin{equation}\label{genClPRODzeros}
t_{n,j,\pm}(p;s) = n^s \big[(2j-1)\pi \pm \arccos \big(-\tfrac{1-p}{p}\big)\big]\, \quad j,n \in\Nset,
\end{equation}
and when $p\in (\tfrac12,1]$ then $\Clg(t)$ changes sign at these zeros.
 So  now we take the logarithm of $|\Clg(t)|$ and study the resulting infinite series of logarithms. 
 This series is the monotone lower limit of a regularized family of series, viz.
\begin{equation}\label{logCloitrePRODgeneralizedREG}
\forall\; t\in\Rset:\quad 
\ln |\Clg(t)| 
= \lim_{\epsilon\downarrow 0} 
{\textstyle\sum\limits_{n\in\Nset}}
\ln \max \left\{\epsilon,\, \left|1- p +p\cos\left(\textstyle n^{-s}_{_{}}t\right)\right|\right\}.
\end{equation}
 For any $\epsilon >0$, the regularized series at the right-hand side can be controlled essentially verbatim to our proof of
the regime $p\in(0,\tfrac12)$. 
 The limit ${\epsilon\downarrow 0}$ is is then established with the help of the integrability of $\ln |t|$ over any 
bounded neighborhood of zero, plus the summability of $1/n^{1+1/s}$ when $1/s>0$.

This ends the outline of our strategy.
\newpage

\noindent
{\textit{Proof of Theorem~\ref{THMgenCLOITREtrend}}:}

	Let $p\in(0,\frac12)$.
 If $t_s^{}>0$ is sufficiently small, then for the \emph{near side} $0\leq t\leq t_s^{}$ 
we have the Maclaurin expansion $\ln\Clg(t)= -\frac{1}{2}p\,\zeta(2s)t^2 +O(t^4)$.
 It follows that $|\ln\Clg(t) + C_{p;s}^{}t^{1/s}|\leq K t^{1/(s+1)}$ for some $K>0$ when $0\leq t \leq t_s^{}$.
 Here and in all estimates below, $K$ is a generic positive constant which may depend on $p,s,\tau,t_s^{}$.

 	As for the \emph{far side} $t\geq t_s$, 
the first sum at r.h.s.(\ref{eq:SdefSsplitSplusONEroot}) is estimated by \vspace{-10pt}
\begin{equation}
\Big|{\textstyle\sum\limits_{n=1}^{{N_s(t/\tau)}}}
\ln\left[1- p +p\cos\left(\textstyle n^{-s}_{_{}}t\right)\right]
\Big|
\; \leq \; \label{eq:firstSsumEST}
	|\ln(1-2p)| \lceil{(st/\tau)^{1/(s+1)}}\rceil
	\; \leq \; K {t^{1/(s+1)}},
\end{equation}

\vspace{-10pt}\noindent
where we used the triangle inequality and
$$
\big|\ln\left[1- p +p\cos\left(\textstyle n^{-s}_{_{}}t\right)\right]\big|\leq |\ln(1-2p)|.
$$ 
	For the second sum at r.h.s.(\ref{eq:SdefSsplitSplusONEroot}) we find, for some $\nu_n\in[n,n+1]$,
\begin{eqnarray}
\textstyle\sum\limits_{n={N_s(t/\tau)+1}}^{\infty}\!\!\!\!\!\!\!\!
\ln\left[1- p +p\cos\left(\textstyle n^{-s}_{_{}}t\right)\right]
-\!\!\displaystyle\int_{N_s(t/\tau)+1}^\infty\!\!\!\!\!\!\!\!\ln\left[1- p +p\cos\left(\textstyle \nu^{-s}_{_{}}t\right)\right] d\nu 
&\!\!\!=\!\!\!& \label{eq:secondSsumESTa}\qquad \\
	 \textstyle\sum\limits_{n={N_s(t/\tau)+1}}^{\infty}\!\!
	\Big(\ln\left[1- p +p\cos\left(\textstyle n^{-s}_{_{}}t\right)\right] - \!\!
\displaystyle\int_n^{n+1}\!\!\!\!\!\!\!\! \ln\left[1- p +p\cos\left(\textstyle \nu^{-s}_{_{}}t\right)\right] d\nu \Big) 
&\!\!\!=\!\!\!& \label{eq:secondSsumESTb}\qquad \\
	\textstyle\sum\limits_{n={N_s(t/\tau)+1}}^{\infty}
\Big(\ln\left[1- p +p\cos\left(\textstyle n^{-s}_{_{}}t\right)\right] -  
        \ln\left[1- p +p\cos\left(\textstyle \nu^{-s}_{n}t\right)\right]\Big)
&\!\!\!=\!\!\!& \label{eq:secondSsumESTc}\qquad \\
	{\textstyle\sum\limits_{n={N_s(t/\tau)+1}}^{\infty}}
\displaystyle \int_{t/\nu_n^{s}}^{t/n^{s}} \frac{p\sin(\xi)}{1-p+p\cos(\xi)}{\rm{d}}\xi;
\end{eqnarray}
here, (\ref{eq:secondSsumESTa}) is obviously true, whereas
(\ref{eq:secondSsumESTb}) expresses the mean value theorem for some $\nu_n\in [n,n+1]$, and 
(\ref{eq:secondSsumESTc}) holds by the fundamental theorem of calculus.
 Now taking absolute values, we estimate
\begin{eqnarray}
	\Big|{\textstyle\sum\limits_{n={N_s(t/\tau)+1}}^{\infty}}
\displaystyle \int_{t/\nu_n^{s}}^{t/n^{s}} \frac{p\sin(\xi)}{1-p+p\cos(\xi)}{\rm{d}}\xi\Big|
&\!\!\!\leq\!\!\!& \label{eq:secondSsumESTd}\qquad \\
	\frac{p}{1-2p}{\textstyle\sum\limits_{n={N_s(t/\tau)+1}}^{\infty}}
\displaystyle	\int_{t/\nu_n^{s}}^{t/n^{s}}	\big|  \sin\xi \big| {\rm{d}}\xi 
&\!\!\!\leq\!\!\!& \label{eq:secondSsumESTe}\qquad \\
	\frac{p}{1-2p}\textstyle\sum\limits_{n={N_s(t/\tau)+1}}^{\infty}
	t\big(\frac{1}{n^{s}}- \frac{1}{\nu_n^{s}}\big) 
&\!\!\!\leq\!\!\!& \label{eq:secondSsumESTf}\qquad \\
\frac{p}{1-2p}	\textstyle\sum\limits_{n={N_s(t/\tau)+1}}^{\infty}
	t\big(\frac{1}{n^{s}}- \frac{1}{(n+1)^{s}}\big) 
&\!\!\!=\!\!\!& \label{eq:secondSsumESTg}\qquad \\
	\frac{p}{1-2p} t\lceil{(st/\tau)^{1/(s+1)}}+1\rceil^{-s}
&\!\!\!\leq\!\!\!& \label{eq:secondSsumESTh}\qquad \vspace{-10pt}\\
 K t^{\frac{1}{s+1}} ;&&\notag
\end{eqnarray}

\vspace{-5pt}\noindent
inequality (\ref{eq:secondSsumESTd}) holds by the triangle inequality in concert with $\cos\xi\geq -1$,
(\ref{eq:secondSsumESTe}) holds since $|\sin \xi|\leq 1$, 
followed by elementary integration, while (\ref{eq:secondSsumESTf}) is due to the monotonic decrease 
of $\nu\mapsto \nu^{-s}$ for $s>1/2$, with  $\nu_n\in [n,n+1]$;
equality (\ref{eq:secondSsumESTg}) holds because the sum at l.h.s.(\ref{eq:secondSsumESTg}) is telescoping; 
lastly, inequality (\ref{eq:secondSsumESTh}) holds because $\ceil{x}$ differs from $x$ by at most 1, and
for large $x$ the $+1$ in its argument becomes negligible.

	For the integral in (\ref{eq:secondSsumESTa}) the variable substitution $\nu^{-s}t=\xi$ yields
\begin{eqnarray}
t^{1/s}	{\textstyle\frac{1}{s}}
\displaystyle\int_0^{t/(N_s(t/\tau)+1)^s}\hspace{-56pt} \ln\left[1- p +p\cos \xi \right]\frac{{\rm{d}}\xi}{\xi^{1+1/s}}  
= \label{eq:SdefSasympINTrewrite} 
t^{1/s}\displaystyle\Biggl[-C_{p;s}^{} - {\textstyle\frac{1}{s}}
	\int_{t/(N_s(t/\tau)+1)^s}^\infty \hspace{-50pt}
 \ln\left[1- p +p\cos \xi \right]\frac{{\rm{d}}\xi}{\xi^{1+1/s}}  \Biggr]\!. 
\end{eqnarray}
Using again the estimate
$| \ln\left[1- p +p\cos \xi \right]|\leq -\ln(1-2p)$, we find (for $t\geq 1$):
\begin{eqnarray}
	t^{1/s}
\Biggr|\int_{t/(N_s(t/\tau)+1)^s}^\infty \hspace{-30pt} \ln\left[1- p +p\cos \xi \right]\frac{1}{s\xi^{1+1/s}} {\rm{d}}\xi\Biggr|
&\!\!\! \leq  &\!\!\!  
\label{eq:SdefSasympINTrewriteESTa}  
	|\ln(1-2p)| \lceil{(st/\tau)^{1/(s+1)}+1}\rceil\  \\ 
& \!\!\!\leq \!\!\! & K t^{1/(s+1)}.
\end{eqnarray}
 This completes the proof of (\ref{genClTHM}) with $|\ln F_{p;s}^{}(|t|)| \leq K_{p;s}^{} |t|^{1/(s+1)}$.

 It remains to prove (\ref{genCLOITREcoeffC}), (\ref{genCLOITREcoeffA}), (\ref{genCLOITREcoeffB}).  
 Integration by parts yields
\begin{equation}
	 C_{p;s}^{} \equiv -{\textstyle\frac1s}\int_0^\infty \ln\left[1- p +p\cos \xi \right]\frac{1}{\xi^{1+1/s}}  {\rm{d}}\xi
= \label{eq:alphaSintPARTIAL}
	 \int_0^\infty \frac{p \sin\xi}{1- p +p\cos \xi}\frac{1}{\xi^{1/s}}  {\rm{d}}\xi,
\end{equation}
where the integral at r.h.s.(\ref{eq:alphaSintPARTIAL}) converges absolutely when $s\in(1/2,1)$, but only conditionally when $s\geq 1$.
 With the help of the geometric series r.h.s.(\ref{eq:alphaSintPARTIAL}) becomes
\begin{equation}
\frac{p}{1- p}  \int_0^\infty \frac{ \sin\xi}{1 + \frac{p}{1- p}\cos \xi}\frac{1}{\xi^{1/s}}  {\rm{d}}\xi 
= \label{Cpq}
\sum_{n=0}^\infty(-1)^n\left(\frac{p}{1-p}\right)^{n+1}\!\! \int_0^\infty {\sin\xi}\cos^n \xi \frac{1}{\xi^{1/s}}  {\rm{d}}\xi;
\end{equation}
the exchange of summation and integration is justified for $s\in (1/2,1)$ by Fubini's theorem, and for $s\geq 1$ by a more
careful limiting argument involving the definition of the conditional convergent integrals as limit $R\to\infty$ of integrals from $0$
to $R$.
 Repeatedly using the trigonometric identity $2\sin(\alpha)\cos(\beta) = \sin(\alpha+\beta)-\sin(\alpha-\beta)$, eventually
followed by a simple rescaling of the integration variable, now yields
\begin{equation}\label{trigINTexpand}
\!\!\!\!\int_0^\infty\!\!\!\!\! {\sin\xi}\cos^n \xi \frac{1}{\xi^{1/s}}  {\rm{d}}\xi = 
\int_0^\infty \!\!\!\!\!{\sin\xi} \frac{1}{\xi^{1/s}}  {\rm{d}}\xi \;
\frac{1}{2^n}\!\!{\textstyle\sum\limits_{k=0}^{\ceil{\frac{n-1}{2}}}}\!\!
\left[\begin{pmatrix} n\cr k \end{pmatrix} -\begin{pmatrix} n\cr k-1 \end{pmatrix}\right]{({1+n -2k})^{\frac1s-1}},
\end{equation}
where it is understood that when $k=0$ one has ${\big(\; {}^{\;n}{}_{\hspace{-10pt}-1}\big)}=0$.
 The integral at r.h.s.(\ref{trigINTexpand}) is $A_{s}^{}$ given in  (\ref{genCLOITREcoeffA}).  
 Inserting (\ref{trigINTexpand}) into (\ref{Cpq}) and pulling out $A_{s}^{}$ yields r.h.s.(\ref{genCLOITREcoeffC}) with
\begin{equation}\label{genCLOITREcoeffBagain}
B_{p;s}^{} 
:= 
\sum_{n=0}^\infty(-1)^n\left(\frac{p}{1-p}\right)^{n+1}\frac{1}{2^n}{\textstyle\sum\limits_{k=0}^{\ceil{\frac{n-1}{2}}}}
\left[\begin{pmatrix} n\cr k \end{pmatrix} -\begin{pmatrix} n\cr k-1 \end{pmatrix}\right]{({1+n -2k})^{\frac1s-1}},
\end{equation}
and a simple manipulation of r.h.s.(\ref{genCLOITREcoeffBagain}) now yields (\ref{genCLOITREcoeffB}).
 This proves the part of Theorem~\ref{THMgenCLOITREtrend} with $p\in(0,\frac12)$.
 
 	Now let $p\in[\frac12,1]$.
        With minor and obvious modifications of our proof for the regime $p\in(0,\frac12)$ one finds that
for $\epsilon>0$ there are $G_{p;s}^{(\epsilon)}(|t|)>0$, $K_{p;s}^{(\epsilon)}>0$, such that
\begin{equation}\label{genClTHMeps}
\forall\,t\in\Rset:\quad 
\prod_{n\in\Nset} \max\left\{\epsilon,\, \left|1-p+p\cos\left({n^{-s}}{t}\right)\right| \right\}= 
\exp\left({- C_{p;s}^{(\epsilon)} \,|t|^{1/s}}\right) G_{p;s}^{(\epsilon)}(|t|),
\end{equation}
with $|\ln G_{p;s}^{(\epsilon)}(|t|)| \leq K_{p;s}^{(\epsilon)} |t|^{1/(s+1)}$, and where
\begin{equation}\label{genCLOITREcoeffCeps}
C_{p;s}^{(\epsilon)} 
:= -\frac1s\int_0^\infty \ln \max\{\epsilon,\,|{1-p+p\cos\xi}|\} \frac{1}{\xi^{1+1/s}}{\rm{d}}\xi.
\end{equation}

 Clearly, $\forall\,t\in\Rset:\;\lim_{\epsilon\downarrow 0}\mbox{l.h.s.}(\ref{genClTHMeps})=
\prod_{n\in\Nset} \left|1-p+p\cos\left({n^{-s}}{t}\right)\right| = |{\mathrm{Cl}}_{p;s}^{}(t)|$. 

 Next, 
$\lim_{\epsilon\downarrow 0}C_{p;s}^{(\epsilon)} =-\frac1s\int_0^\infty \ln |{1-p+p\cos\xi}| \frac{1}{\xi^{1+1/s}}{\rm{d}}\xi<\infty$
as a doubly improper Riemann integral.
 Indeed, the singularities, one at each zero of $1-p+p\cos\left(\xi\right)$, are all improper
Riemann integrable because $-\int_0^\delta \ln \xi \rm{d}\xi = -\delta\ln\delta +\delta \downarrow 0$ as $\delta\downarrow 0$.
 There are countably many singularities, located at 
\begin{equation}\label{genClPRODzerosXI}
\xi_{j,\pm} = (2j-1)\pi \pm \arccos \big(-\tfrac{1-p}{p}\big)\, \quad j\in\Nset,
\end{equation}
and so the absolute contribution from a $\delta$-neighborhood of the singularity at $\xi_{j,\pm}$, 
with $\delta\downarrow 0$ when $\epsilon\downarrow 0$, can be dominated by $c(-\delta\ln\delta +\delta)/ \xi_{j,\pm}^{1+1/s}$
for some positive constant $c$ which is independent of $j$ and the $\pm$ index.
 Since $1/j^{1+1/s}$ is summable over $\Nset$ when $1/s>0$, the absolute difference between $C_{p;s}$ and 
$C_{p;s}^{(\epsilon)}$ is dominated by $\tilde{c} (-\delta\ln\delta +\delta)$, which vanishes as 
$\delta\downarrow 0$ when $\epsilon\downarrow 0$.

 The convergence of l.h.s.(\ref{genClTHMeps}) and of the first factor at r.h.s.(\ref{genClTHMeps}) now imply that also 
 the second factor at r.h.s.(\ref{genClTHMeps}) converges when $\epsilon\downarrow 0$.
 This does not yet establish an upper bound on some $|F_{p;s}^{}(|t|)|$, defined below, as claimed in the theorem;
indeed, the bound $|\ln G_{p;s}^{(\epsilon)}(|t|)| \leq K_{p;s}^{(\epsilon)} |t|^{1/(s+1)}$ and the fact
that $G_{p;s}^{(\epsilon)}(|t|)$ has zeros in the limit $\epsilon\to0$ means that 
the constants $K_{p;s}^{(\epsilon)}\to\infty$ as $\epsilon\to0$. 
 However, our theorem only requires an {upper bound} on $\ln G_{p;s}^{(\epsilon)}(|t|)$,
\emph{not} on $|\ln G_{p;s}^{(\epsilon)}(|t|)|$, as ${\epsilon\to0}$.
 To prove such a bound is straightforward. 
 Since l.h.s.(\ref{genClTHMeps}) only introduces a lower cutoff on $|\Clg(t)|$, at $\epsilon$, it follows by 
inspection of the estimates in the proof of the $p\in(0,\frac12)$ part of our Theorem that 
we do have the upper bound $\ln G_{p;s}^{(\epsilon)}(|t|) \leq K_{p;s}^{} |t|^{1/(s+1)}$ 
\emph{uniformly} in $\epsilon$, and this proves that $G_{p;s}^{}(|t|) \leq \exp(K_{p;s}^{} |t|^{1/(s+1)})$ in the
limit $\epsilon\downarrow 0$.
 Finally, we set $F_{p;s}(|t|) := G_{p;s}^{}(|t|)\; \sign\,{\mathrm{Cl}}_{p;s}^{}(t)$, and
the entirely elementary proof of Theorem~\ref{THMgenCLOITREtrend} is complete. \hfill QED

\newpage

\section{L\'evy Trends and Fluctuations}\label{TandF}

 In this section we display the PDFs for a small selection of random Riemann-$\zeta$ walks $\Omega_{p}^{\zeta}(s)$,
obtained by numerical Fourier transform of their characteristic functions $\Clg(t)$.
 We compare them with the Fourier transform of their trend functions $\exp\left(- C_{p;s}^{} \,|t|^{1/s}\right)$,
which are known as L\'evy-stable distributions with \emph{stability parameter} $\alpha=1/s$, 
\emph{skewness parameter} $\beta=0$, \emph{scale parameter} $c=C_{p;s}^{s}$, and \emph{median} $\mu=0$; see \cite{probBOOK}.
 The comparison will highlight the importance of the fluctuating factors 
$F_{p;s}(|t|)$ 
in the characteristic functions $\Clg(t)$.

 The first figure shows the PDF $f_{\Omega_{p}^\zeta(s)}(\omega)$ 
for Cloitre's parameter values $p=1/3$ and $s=2$, together with the pertinent L\'evy PDF (here $\cC$ and $\cS$ are Fresnel integrals)

\begin{equation}\label{LevyHALFetc}
f_{\Omega_{\hfrac12;0;C^{2};0}^{\mbox{\tiny{\sc{L\'evy}}}}}\!(\omega) =
2\pi u^3 \Bigl(\sin\!\left(\tfrac{\pi}{2}u^2\right)\!\!\left[\tfrac12-\cS\!\left(u\right)\right]
             +\cos\!\left(\tfrac{\pi}{2}u^2\right)\!\! \left[\tfrac12-\cC\!\left(u\right)\right]\Bigr),
\end{equation}

\noindent
where $u = C/\sqrt{2\pi|\omega|}$ and $C=C_{1/3;2}$; cf. the histogram Fig.~4.
\vspace{5pt}

\includegraphics[scale=0.65]{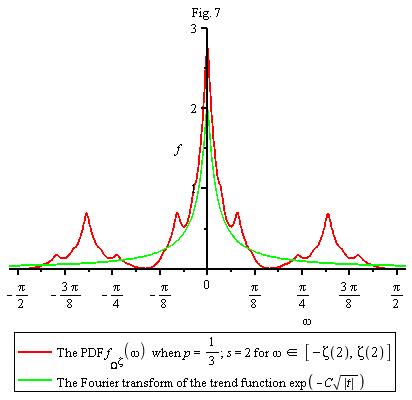}

\vspace{5pt}
\noindent
 Fig.~7  reveals that the stable distribution (\ref{LevyHALFetc})
obtained by Fourier transform of the L\'evy trend factor $\exp({- C\surd|t|})$, 
which captures the ``large scale'' behavior of $\Cl(t)$ asymptotically exactly but misses all of 
its ``small scale'' details (recall Fig.~1 and Fig.~2), only very crudely resembles the distribution 
obtained by the Fourier transform of $\Cl(t)$.
 Also, we recall that the random variable $\Omega_{1/3}^\zeta(2)$ takes its values in the interval $[-\zeta(2),\zeta(2)]$, so
$f_{\Omega_{1/3}^\zeta(2)}(\omega)$ vanishes identically outside this interval. 
 By contrast, L\'evy-stable PDF are ``heavy-tailed'' (except when $\alpha=2$, i.e. $s=1/2$, which is excluded here); in particular,
it follows from (\ref{LevyHALFetc}) (see also \cite{probBOOK}) that
\begin{equation}\label{LevyHALFetcASYMP}
f_{\Omega_{\frac12;0;C;0}^{\mbox{\tiny{\sc{L\'evy}}}}}\!(\omega) 
\sim 
\tfrac12 C_{\frac13;2}^{2} {\sqrt{\pi}}|\omega|^{-3/2}
\quad (\omega\to\infty).
\end{equation}

 Next we turn to the borderline case $s=1$, which is particularly interesting. 
 When $p\neq 1$ this random walk is a  generalization of the harmonic random walk ($p=1$) studied by Kac \cite{Kac},
 Morrison \cite{Morrison}, and Schmuland \cite{Schmuland}.
 Furthermore,  
the ``trend factor'' of the characteristic function for ${\Omega_{p}^\zeta(1)}^{}$ 
becomes $e^{-C_{p;1}^{}|t|}$: the characteristic function of a Cauchy random variable with ``theoretical spread'' 
$C_{p;1}^{}$ (which is explicitly computable; see below).
 The next Figure displays the PDF $f_{\Omega_{1/3}^\zeta(1)}(\omega)$ for the harmonic random walk with $p=1/3$,
together with the Cauchy distribution of theoretical spread $C_{1/3;1}^{}$ about~$0$; cf. the histogram in Fig.~5.

 \smallskip

\includegraphics[scale=0.55]{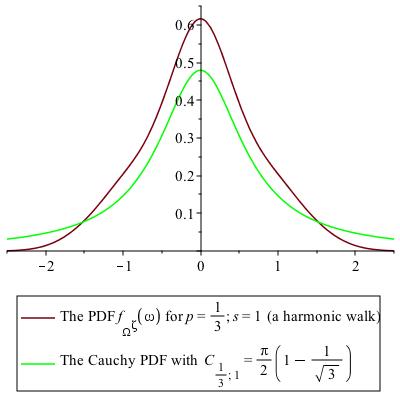}

 The discrepancy between the PDF $f_{\Omega_{1/3}^\zeta(1)}(\omega)$ for the harmonic random walk with $p=\frac13$
and the Cauchy distribution of theoretical spread $C_{1/3;1}^{}$ about~$0$ 
visible in Fig.~8 is not quite as flagrant as the corresponding discrepancy in Fig.~7.
 Not so outside the shown interval, though: the Cauchy distribution is heavy-tailed, while  
$f_{\Omega_{1/3}^\zeta(1)}(\omega)$, because it is a Schwartz function, has moments of all order.
 This can also be shown by adaptation of the estimates for $f_{\Omega_{1/2}^\zeta(1)}(\omega)$ given by
Schmuland \cite{Schmuland}.

 We also vindicate our claim that one can compute $C_{p;1}^{}$ explicitly.
 First of all, 
\begin{equation}\label{genCLOITREcoeffWHENsISoneBa}
\sum_{k=0}^{\ceil{\hfrac{(n-1)}{2}}}
\begin{pmatrix} n\cr k \end{pmatrix} \frac{1+n-2k}{{1+n-k}\;}
=
\begin{pmatrix} n\cr \floor{\hfrac{n}{2}}\end{pmatrix},
\end{equation}
which is A001405 in Sloane's OEIS. 
 Now $\frac{1}{2^0}\genfrac{(}{)}{0pt}{}{0}{\floor{\hfrac{0}{2}}} = 1$ while 
$\frac{1}{2^n}\genfrac{(}{)}{0pt}{}{n}{\floor{\hfrac{n}{2}}}= \frac{1}{2^{n-1}}\genfrac{(}{)}{0pt}{}{n-1}{\floor{\hfrac{n-1}{2}}}$ 
when $n=2m$ with $m\in\Nset$, and using that 
$\sum_{m=0}^\infty \! \frac{1}{2^{2m}} \genfrac{(}{)}{0pt}{}{2m}{m} x^{2m} = \frac{1}{\sqrt{1-x^2}}$
we compute  

\vspace{-15pt}
\begin{eqnarray}\label{genCLOITREcoeffWHENsISoneBb}
B_{p;1}^{} =\! \sum_{n=0}^\infty\left(-1\right)^n\!\left(\!\frac{p}{1-p}\!\right)^{\!\!n+1}\frac{1}{2^n}
\begin{pmatrix} n\cr \floor{\hfrac{n}{2}}\end{pmatrix}
= 
1 - \sqrt{1-2p}\quad \mbox{for}\quad  p\in(0,\tfrac12);
\end{eqnarray}

\vspace{-10pt}\noindent
so with $A_1^{}=\frac{\pi}{2}$ we obtain $C_{p;1}^{} = A_1^{} B_{p;1}^{}$ in closed form, displayed in Fig.~9.
 Note that its $p$-derivative blows up as $p\nearrow\frac12$.
\smallskip

\includegraphics[scale=0.5]{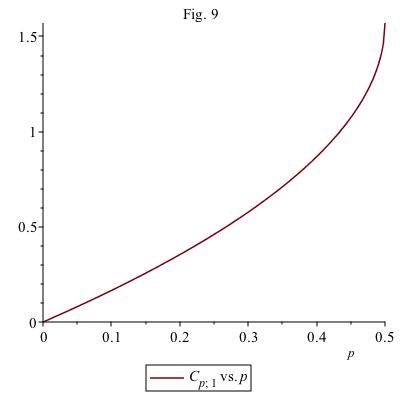}

\section{Open Problems}\label{FIN}

 The following problems seem to be particularly worthy of further pursuit.\vspace{-15pt}
 
\subsection{Why L\'evy trends?}\label{LevySECRETS}\vspace{-5pt}

 What is the probabilistic reason for the occurrence of the symmetric L\'evy $\frac1s$-stable 
distributions associated with the trend factors?
 We recall that $X$ is a \emph{L\'evy-stable} random variable if and only if
$X= c_1 X_1 +c_2X_2$, where $X_1$ and $X_2$ are i.i.d. copies of $X$ and $c_1$ and $c_2$ are suitable
positive constants; see also \cite{GaroniFrankel}.
 Where is this ``L\'evy stability'' hiding in the random Riemann-$\zeta$ walks?\vspace{-10pt}

\subsection{Does the singularity at $p=\tfrac12$ have statistical significance?}\label{singularSECRETS}\vspace{-5pt}
 The derivative singularity of $p\mapsto C_{p;s}$ at $p=\frac12$ is inherited from the derivative singularity of the absolute value function.
 Is this a consequence of our method of representing $\Clg(t)$, or does this have some statistical physics meaning for the family of random walks?
 Something akin to a ``percolation threshold''? 
 In the random walks with $p<\frac12$ one more often stays put than moving to another position, while for $p>\frac12$ the opposite is true. 
 Does this entail a singular change in the statistical random walk behavior, or is this only a peculiar singularity in the trend function?
 
\subsection{Are there ``perfectly typical'' random Riemann-$\zeta$ walks?}\label{typicalSECRETS}\vspace{-5pt}
 If the intersection of all typical subsets of the set of random Riemann-$\zeta$ walks 
for given $p\in(0,1]$ and $s>0$ is non-empty, then the answer is ``Yes!'' --- 
in that case it would be very interesting to exhibit a perfectly typical walk explicitly, if at all possible.
 It is also conceivable that the intersection set is empty.\vspace{-10pt}

\subsection{Complex random Riemann-$\zeta$ walks}\label{complexSECRETS}\vspace{-5pt}
 What happens if one extends $\Omega_{p}^\zeta(s)$ to complex $s$?
 The Riemann hypothesis implies for $\zeta(s)$ itself that its extremal walks with Im$(s)\neq 0$
converge to the origin if and only if $\Re{e}(s)=\hfrac12$ and Im$(s)$ is the imaginary part of a nontrivial zero of $\zeta(s)$. 
 Does $\Re{e}(s)=\hfrac12$ play a special role also for the random Riemann-$\zeta$ walks?

 \vfill

\section*{Acknowledgement}\vspace{-10pt}
We truly thank: Benoit Cloitre for posing his problem; Alex Kontorovich for his enlightening explanations of
the Riemann hypothesis; 
Norm Frankel and Larry Glasser for their interest in and helpful feedback on $C_{p;s}^{}$;
Neil Sloane for OEIS and for his comments; Doron Zeilberger for noting that the combinatorics 
in our evaluation of $C_{p;s}^{}$ produces Catalan numbers. 
 We also thank the referees for constructive comments.
 Some symbolic manipulations were obtained with MAPLE, as were the figures.

\vfill
 \newpage

\section*{Appendix on Power Walks}

 If instead of a step size which decreases by the power law $n\mapsto n^{-s}$ one uses an exponentially
decreasing step size $n\mapsto s^{-n}$ with $s>1$, the outcome is a ``random geometric series'' 
(a sum over powers of $1/s$ with random coefficients $R_p^{}(n)\in\{-1,0,1\}$),

\vspace{-15pt}
\begin{equation}\label{RpowSERIES}
\Omega_{p}^{\mbox{\tiny{pow}}}(s) := {\sum_{n\in\Nset}^{}} R_p^{}(n)\frac{1}{s^n}, \quad s >1,\quad p\in(0,1];
\end{equation}

\vspace{-5pt}
\noindent
the pertinent walks are called ``geometric walks.''
 With more general random coefficients one simply speaks of ``random power series'' and their ``power walks.'' 

 All these random variables $\Omega_{p}^{\mbox{\tiny{pow}}}(s)$ have characteristic functions with infinite 
trigonometric product representations obtainable from our (\ref{charFUNC}) by 
replacing ${}^\zeta\to {}^{\mbox{\tiny{pow}}}$ and $n^{-s}\to s^{-n}$.
 Some of these can be evaluated  in terms of elementary functions.
 We register a few special cases, beginning with three geometric walks and ending with a countable family of more general (but simple)
power walks.
\smallskip

\noindent
(i) Setting $p=1$ and $s=2$ gives the chacteristic function (see formula (1) of \cite{Morrison}) 

\vspace{-15pt}
\begin{equation}\label{EULERprod}
\Phi^{}_{\Omega_{1}^{\mbox{\tiny{pow}}}(2)}(t)
 = \prod_{n\in\Nset} \cos\left(\frac{t}{2^n}\right) \equiv \frac{\sin t}{t},
\end{equation}

\vspace{-5pt}
\noindent
an infinite product\footnote{By substituting $\pi/2$ for $t$ and repeatedly 
using a trigonometric angle-halving identity one arrives at Vi\`ete's infinite product for $2/\pi$, allegedly the 
first infinite product ever proposed.}
representation of the sinc function derived by Euler algebraically by exploiting the trigonometric angle-doubling formulas 
(see \cite{Morrison}).
 Recall that sinc$(t)=\int_{-1}^1\frac12 e^{it\omega}d\omega$ is the (inverse) Fourier transform of the 
PDF $f_{\Omega^{\mbox{\tiny{unif}}}}(\omega)$ of the uniform random variable $\Omega^{\mbox{\tiny{unif}}}$ on $[-1,1]$, i.e.
$f_{\Omega^{\mbox{\tiny{unif}}}}(\omega) = \frac12$ if $\omega \in [-1,1]$, and $f_{\Omega^{\mbox{\tiny{unif}}}}(\omega) =0$ otherwise.
 Indeed, $\Omega_{1}^{\mbox{\tiny{pow}}}(2)$ is a random walk representation of $\Omega^{\mbox{\tiny{unif}}}$ 
equivalent to the binary representation of $[0,1]$: recalling that any real number $x\in[0,1]$ has a 
binary representation\footnote{Those representations are not unique and one needs to consider their equivalence 
classes to identify them uniquely with their real outcome on $[0,1]$, cf. \cite{Kac}.}
$x = 0.b_1b_2b_3...\equiv \sum_{n\in\Nset} b_n /2^n$ with $b_n\in\{0,1\}$, and 
noting that if $x\in[0,1]$ then $\omega:=2x-1\in[-1,1]$, 
it follows that any real number $\omega\in[-1,1]$ has a binary representation
$\omega =  \sum_{n\in\Nset} r_2^{}(n) /2^n$ with $r^{}_2(n)\in\{-1,1\}$.
 It is manifest that any such representation of $\omega$ is an outcome of $\Omega_{1}^{\mbox{\tiny{pow}}}(2)$.
\smallskip

\noindent
(ii) $\Omega_{1}^{\mbox{\tiny{pow}}}(3)$ is the random variable for which the characteristic function

\vspace{-15pt}
\begin{equation}\label{MorrisonCANTORprod}
\Phi^{}_{\Omega_{1}^{\mbox{\tiny{pow}}}(3)}(t)
= \prod_{n\in\Nset} \cos\left(\frac{t}{3^n}\right) =:\Phi^{}_{\Omega^{\mbox{\tiny\sc{Cantor}}}}(t/2)
\end{equation}

\vspace{-5pt}
\noindent
is a trigonometric product discussed in \cite{Morrison}.
 Morrison explains that $\Phi^{}_{\Omega^{\mbox{\tiny\sc{Cantor}}}}(t)$ is the characteristic function of a
random variable ${\Omega^{\mbox{\tiny\sc{Cantor}}}}$ that is uniformly distributed over the
Cantor set constructed from $[-1,1]$ by removing middle thirds ad infinitum.
 For uniform distributions on other Cantor sets, see \cite{DFT}.

 We remark that this is a nice example of a random walk whose endpoints are distributed by a singular distribution, 
in the sense that the Cantor set obtained from $[-1,1]$ has Lebesgue measure $0$. 
 As pointed out to us by an anonymous referee, the distribution of $\Omega_{1}^{\mbox{\tiny{pow}}}(s)$ 
is singular and concentrated on some Cantor set \emph{for all} $s>2$, while
for $1 < s < 2$ the story is more complicated: 
  Solomyak \cite{Solomyak} proved that the distribution of $\Omega_{1}^{\mbox{\tiny{pow}}}(s)$ 
is absolutely continuous (i.e., it is equivalent to a PDF, an integrable function)
\emph{for almost every} $s \in (1, 2)$; see also \cite{PerezSolomyak}.
 However, the distribution of $\Omega_{1}^{\mbox{\tiny{pow}}}(s)$ is not absolutely continuous for all $s\in(1,2)$ ---
in 1939 Erd\H{o}s found values of $s\in(1,2)$ for which the distribution is singular; these are still the only ones known. 
 See \cite{PerezSchlagSolomyak} for further reading.
\smallskip

\noindent
(iii)  Setting $p=\frac23$ and $s=3$ yields 

\vspace{-15pt}
\begin{equation}\label{MorrisonPROD}
\Phi^{}_{\Omega_{2/3}^{\mbox{\tiny{pow}}}(3)}(t) = 
\prod_{n\in\Nset}\left[\frac13 +\frac23\cos\left(\frac{t}{3^{n}}\right)\right] \equiv
\frac{\sin (t/2)}{t/2} ,
\end{equation}

\vspace{-5pt}
\noindent
which becomes formula (9) of \cite{Morrison} under the rescaling $t\mapsto 2t$
(see also Exercise 3 on page 11 of \cite{Kac}.)
 Recalling our discussion of example (i), we conclude that  (\ref{MorrisonPROD}) is the characteristic function of the
uniform random variable on the interval $[-\frac12,\frac12]$, expressed as a random walk equivalent to the ternary 
representation of the real numbers in $[0,1]$, shifted to the left by $-\frac12$.
 
\noindent
(iv) The sinc representations (\ref{EULERprod}) and (\ref{MorrisonPROD}) (after rescaling $t\mapsto 2t$) 
are merely the first two members of a countable family of infinite trigonometric product representations of $\sin t / t$
derived by Kent Morrison \cite{Morrison}, and given by
\begin{equation}\label{MorrisonPRODgeneral}
 \frac{\sin t}{t} =  \prod_{n\in\Nset}\sum_{m=1-s}^{s-1}\frac{1-(-1)^{s+m}\hspace{-15pt}}{2s}\hspace{15pt} 
\cos\left(\frac{m}{s^n}t\right),\; 1<s\in\Nset;
\end{equation}
$s$ even in (\ref{MorrisonPRODgeneral}) is formula (12) in \cite{Morrison}, 
$s$ odd  in (\ref{MorrisonPRODgeneral}) is formula (13) in \cite{Morrison}.
 These representations of the characteristic function of the uniform random variable over $[-1,1]$ are obtained by
considering random walks that enter with equal likelihood into any one of $s$ branches which ``$s$-furkate'' off of every vertex
of a symmetric tree centered at 0, equivalent to the usual ``$s$-ary'' representation of the
real numbers in $[0,1]$ (shifted to the left by $-\frac12$ and scaled up by a factor 2).
 When $s>3$ these are no longer random geometric series, but still simple random power series.

\noindent


\vfill\vfill
\end{document}